\documentclass{article}
\usepackage{amssymb,amsmath,amsthm,graphicx,cite}

\def\qed{\hfill {\hbox{${\vcenter{\vbox{               
   \hrule height 0.4pt\hbox{\vrule width 0.4pt height 6pt
   \kern5pt\vrule width 0.4pt}\hrule height 0.4pt}}}$}}}

\def\utr{\, \underline{\triangleright}\, }
\def\otr{\, \overline{\triangleright}\, }

\newtheorem{theorem}{Theorem}

\newtheorem{corollary}[theorem]{Corollary}

\newtheorem{observation}{Observation}
\theoremstyle{definition}
\newtheorem{example}{Example}
\newtheorem{definition}{Definition}

\date{}

\title{\Large \textbf{Biquandle Virtual Brackets}}

\author{Sam Nelson\footnote{Email: Sam.Nelson@cmc.edu. Partially supported by Simons Foundation collaboration grant 316709}\and
Kanako Oshiro\footnote{Email: oshirok@sophia.ac.jp.
Partially supported by JSPS KAKENHI Grant Number 16K17600. } \and
Ayaka Shimizu\footnote{Email: shimizu@nat.gunma-ct.ac.jp.
Partially supported by Grant for Basic Science Research Projects from The Sumitomo Foundation (160154).} \and
Yoshiro Yaguchi\footnote{Email: yaguchi-y@nat.gunma-ct.ac.jp}}

\begin{document}
\maketitle

\begin{abstract}
We introduce an infinite family of quantum enhancements of the biquandle 
counting invariant we call \textit{biquandle virtual brackets}. Defined in 
terms of skein invariants of biquandle colored oriented knot and link diagrams 
with values in a commutative ring $R$ using virtual crossings as smoothings, 
these invariants take the form of multisets of elements of $R$ and can be 
written in a ``polynomial'' form for convenience. The family of invariants 
defined herein includes as special cases all quandle and biquandle 2-cocycle 
invariants, all classical skein invariants (Alexander--Conway, Jones, HOMFLYPT 
and Kauffman polynomials) and all biquandle bracket invariants defined in 
\cite{NOR} as well as new invariants defined using virtual crossings in a 
fundamental way, without an obvious purely classical definition.
\end{abstract}

\parbox{6in} {\textsc{Keywords:} Quantum enhancements, biquandles, biquandle
counting invariants, virtual knots and links

\smallskip

}

\section{\large\textbf{Introduction}}\label{I}

In \cite{NR} the notion was introduced of \textit{quantum enhancements} of 
the biquandle counting invariant by means of quantum invariants of 
biquandle-colored knots and links. More precisely, if $\beta$ is a quantum 
invariant of biquandle colored knots and links, then the multiset of 
$\beta$-values over the set of biquandle colorings of a knot or link $K$
is an invariant whose cardinality recovers the biquandle counting invariant.
In \cite{NOR}, a family of such quantum enhancements including both 
classical skein invariants and biquandle cocycle invariants as special cases,
known as \textit{biquandle brackets}, was introduced. These invariants can be 
understood as skein invariants of biquandle colored oriented knots and links
with skein coefficients depending on the biquandle colors of the semiarcs 
involved in the crossing being smoothed. The fact that a biquandle coloring 
is ``broken'' by smoothings means it is simpler to think of the invariant in 
terms of the state-sum definition in which all smoothings are done at once 
rather than smoothing crossings one at a time. In \cite{IM}, picture-valued 
biquandle brackets were introduced, and in \cite{NO} \textit{trace diagrams}
are used to compute biquandle brackets via a recursive expansion instead of
via the state-sum definition.

In this paper we generalize the biquandle bracket idea to use a different
set of skein relations in which a virtual crossing is considered as a type 
of a smoothing rather than a kind of crossing. We use the 
term ``biquandle virtual bracket'' to distinguish this 
case from ``virtual biquandle brackets'' which are brackets with coefficients 
in a virtual biquandle. This new infinite family of invariants of oriented 
classical and virtual knots and links contains classical biquandle brackets
(and hence classical skein invariants and biquandle 2-cocycle invariants)
as special cases and uses virtual knot theory in a fundamental way, joining 
the finite type invariants of \cite{GPV} as a family of invariants of 
classical knots without an obvious definition in term of purely classical 
knot and link diagrams.

The paper is organized as follows. In Section \ref{B} we briefly review 
biquandles, biquandle colorings of oriented classical and virtual knots and 
links, and the biquandle counting invariant. In Section \ref{VB} we introduce
biquandle virtual brackets. In Section \ref{AE} we provide some applications 
and examples. In Section \ref{Q} we close with some questions for future 
research.

\section{\large\textbf{Biquandles and Virtual Links}}\label{B}

We begin with a definition (see \cite{EN,FJK,FRS}).

\begin{definition}
A \textit{biquandle} is a set $X$ with operations $\utr,\otr:X\to X$
satisfying for all $x,y,z\in X$
\begin{itemize}
\item[(i)] $x\utr x=x\otr x$,
\item[(ii)] the maps $\alpha_x,\beta_x:X\to X$ and $S:X\times X\to X\times X$
given by 
\[\alpha_x(y)=y\otr x,\quad \beta_x(y)=y\utr x\quad \mathrm{and}\quad
S(x,y)=(y\otr x, x\utr y)\]
are invertible, and
\item[(iii)] we have the \textit{exchange laws}
\[\begin{array}{rcl}
(x\utr y)\utr(z\utr y) & = & (x\utr z)\utr(y\otr z) \\
(x\utr y)\otr(z\utr y) & = & (x\otr z)\utr(y\otr z) \\
(x\otr y)\otr(z\otr y) & = & (x\otr z)\otr(y\utr z).
\end{array}\]
\end{itemize}
It is sometimes convenient to
write $x\utr y$ as $x^y$ and $x\otr y$ as $x_y$ for the sake of space.
\end{definition}

\begin{definition}
A map $f:X\to Y$ between biquandles is a \textit{biquandle homomorphism}
if 
\[f(x\utr y)=f(x)\utr f(y)\quad \mathrm{and}\quad f(x\otr y)=f(x)\otr f(y)\] 
for all
$x,y\in X$. A subset $S\subset X$ of a biquandle $X$ is a \textit{subbiquandle}
if $S$ is closed under the biquandle operations of $X$, including the 
inverses of the maps in axiom (ii).
\end{definition}

\begin{example}
A biquandle $X$ in which $x\otr y=x$ for all $y\in X$ is a \textit{quandle}. 
Examples of quandles include
\begin{itemize}
\item \textit{$n$-fold conjugation quandles:} Every group $G$ (or union of 
conjugacy classes in a group) is a quandle under $x\utr y=y^{-n}xy^n$ for 
$n\in\mathbb{Z}$,
\item \textit{Alexander quandles:} Every module $X$ over $\mathbb{Z}[t^{\pm 1}]$
is a quandle under $x\utr y=tx+(1-t)y$,
\item \textit{symplectic quandles:} Every vector space $V$ over a field 
$\mathbb{F}$ of characteristic $\ne 2$ with a symplectic form 
$\langle,\rangle:V\times V\to \mathbb{F}$ is a quandle  under
$x\utr y=x+\langle x,y\rangle y$. 
\end{itemize}
\end{example}

\begin{example}
Examples of non-quandle biquandles include
\begin{itemize}
\item \textit{constant action biquandles:} Every set $X$ is a biquandle 
under $x\otr y=x\utr y=\sigma(x)$ for a bijection $\sigma:X\to X$,
\item \textit{Alexander biquandles:} Every module $X$ over the ring 
$\mathbb{Z}[t^{\pm 1},s^{\pm 1}]$ is a biquandle under the operations
$x\utr y=tx+(s-t)y$ and $x\utr y=sx$.
\item \textit{fundamental biquandle of an oriented link:} The set 
$\mathcal{B}(L)$ of equivalence classes of \textit{biquandle words} in 
generators corresponding to semiarcs in a diagram of an oriented link $L$ 
modulo the congruence generated by the crossing relations
\[\includegraphics{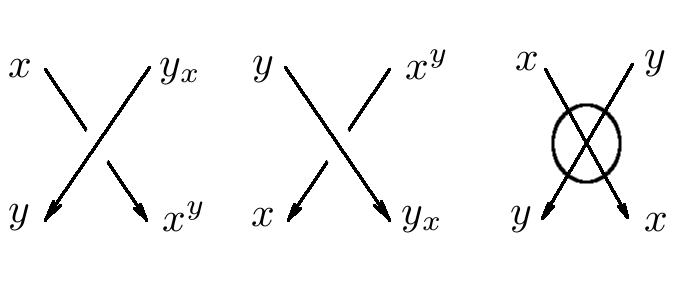}\]
and the biquandle axioms is a strong invariant of virtual links, conjectured 
to be complete up to reflection for virtual knots and known to be complete up 
to reflection for classical knots (see \cite{EN,FJK,HK} for more).
\end{itemize}
\end{example}

Given a finite set $X=\{x_1,\dots,x_n\}$ we can define a biquandle structure
on $X$ by listing operation tables of $\utr$ and $\otr$ such that the 
axioms are satisfied -- the diagonals of the tables must agree, the columns 
must be permutations, etc. We can conveniently encode these tables as an 
$n\times 2n$ block matrix by dropping the ``$x$''s and listing only the 
subscripts, e.g.
\[
\begin{array}{r|rrr}
\utr & x_1 & x_2 & x_3 \\ \hline
x_1 & x_1 & x_3 & x_2 \\
x_2 & x_3 & x_2 & x_1 \\
x_3 & x_2 & x_1 & x_3
\end{array}\quad 
\begin{array}{r|rrr}
\otr & x_1 & x_2 & x_3 \\ \hline
x_1 & x_1 & x_1 & x_1 \\
x_2 & x_2 & x_2 & x_2 \\
x_3 & x_3 & x_3 & x_3
\end{array}
\rightarrow 
\left[
\begin{array}{rrr|rrr}
 1 & 3 & 2 & 1 & 1 & 1 \\
 3 & 2 & 1 & 2 & 2 & 2\\
 2 & 1 & 3 & 3 & 3 & 3
\end{array}
\right]
\]

Next, we recall the definition of virtual knots and links (see \cite{K} for 
more). A \textit{virtual knot diagram} has classical and virtual crossings
as depicted:
\[\includegraphics{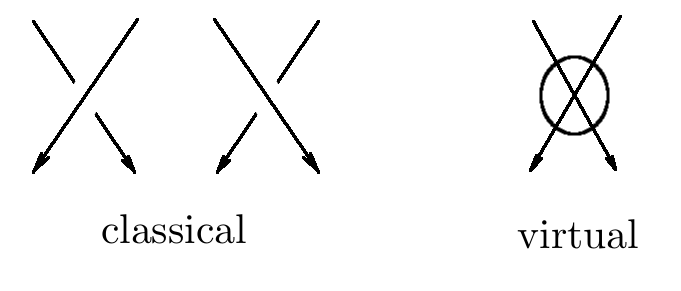}\]
The virtual crossings are interpreted as not really crossings, but  
artifacts of drawing a non-planar knot diagram on planar paper -- we can
regard a virtual crossing as indicating genus in the supporting surface on 
which the knot diagram is drawn. In this way, a virtual crossing may be 
considered a kind of smoothing rather than a kind of crossing.

Two virtual link diagrams are equivalent if they are related by the 
\textit{oriented virtual Reidemeister moves:}
\[\includegraphics{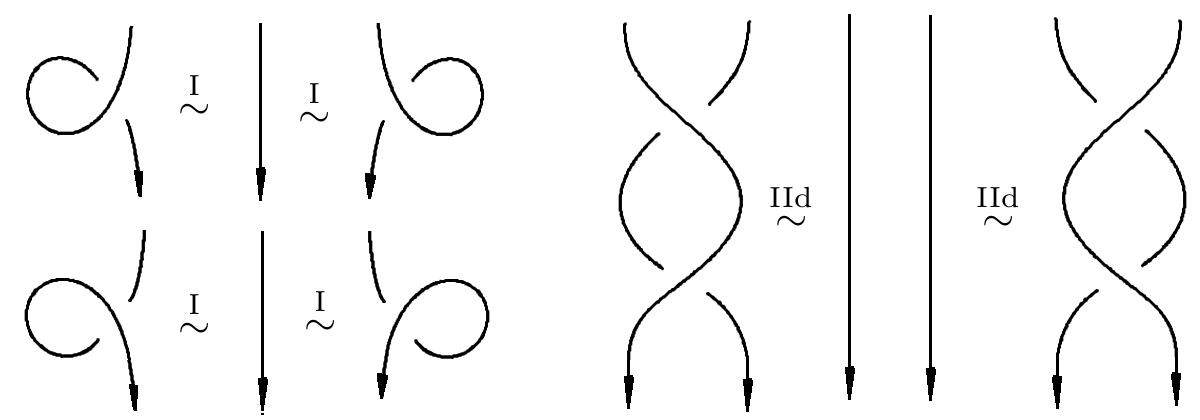}\]
\[\includegraphics{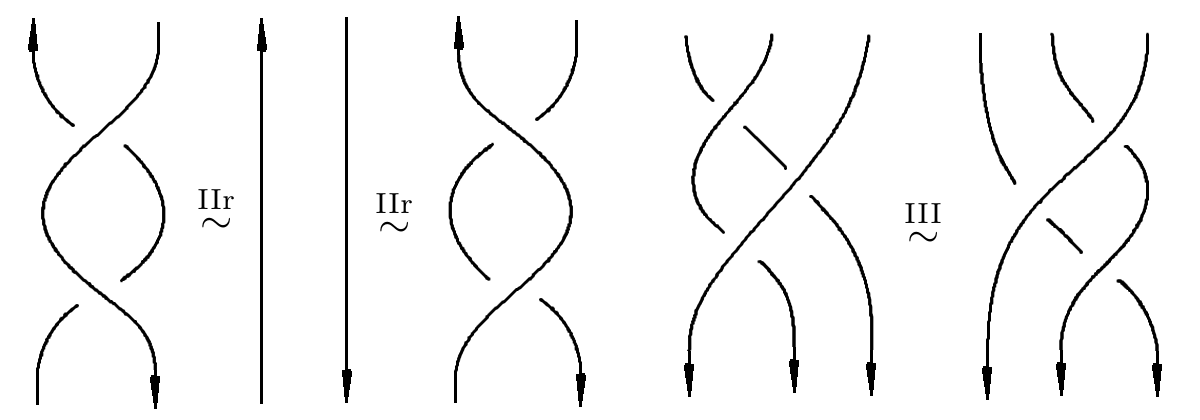}\]
\[\includegraphics{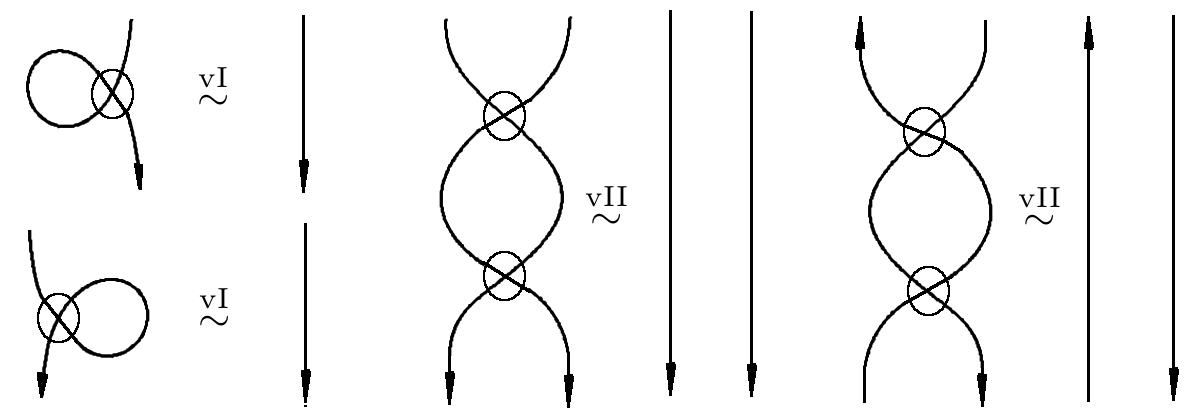}\]
\[\includegraphics{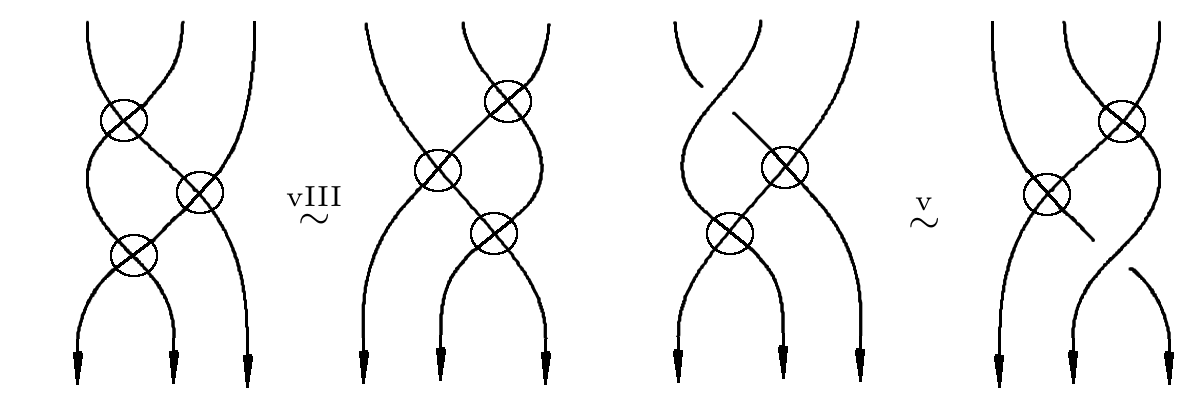}\]
Note that there are more oriented III, vIII and v moves, but in the presence 
of the other pictured moves they can be recovered from the pictured moves. In 
particular, the above constitute one generating set of oriented virtual 
Reidemeister moves. Classical knots and links form a subset of virtual knots 
and links, with two classical knots or link equivalent as virtual knots or 
links if and only if they are equivalent as classical knots or links; see 
\cite{K} for more.

Given a biquandle $X$ and an oriented virtual knot or link $L$ represented by 
a diagram $D$, we can \textit{color} $D$ with $X$ by assigning elements of 
$X$ to the semiarcs in $D$, i.e. the portions of $D$ between crossing points, 
such that we locally have the following pictures:
\[\includegraphics{sn-ko-as-yy-12.png}\]
It is then straightforward to verify that for any $X$-coloring of $D$ before
a Reidemeister move, there is a unique corresponding coloring after the move.
In particular, if $X$ is finite, the number of $X$-colorings of $D$
is finite.
In terms of the fundamental biquandle, a coloring is a homomorphism
$f:\mathcal{B}(L)\to X$, and the set of colorings is 
$\mathrm{Hom}(\mathcal{B}(L),X)$; then the counting invariant is
$\Phi_X^{\mathbb{Z}}(L)=|\mathrm{Hom}(\mathcal{B}(L),X)|$.

Hence, we have

\begin{theorem}
The number of $X$-colorings of a virtual knot or link $K$ by a biquandle $X$,
denoted $\Phi_X^{\mathbb{Z}}(K)$ and known as the \textit{biquandle counting 
invariant}, is an invariant of oriented classical and virtual links. 
\end{theorem}

\begin{example}\label{ex:HL}
The \textit{Hopf link} below has four colorings by the biquandle 
$X=\mathbb{Z}_2$ with $x\utr y=x\otr y=x+1$ as depicted:
\[\includegraphics{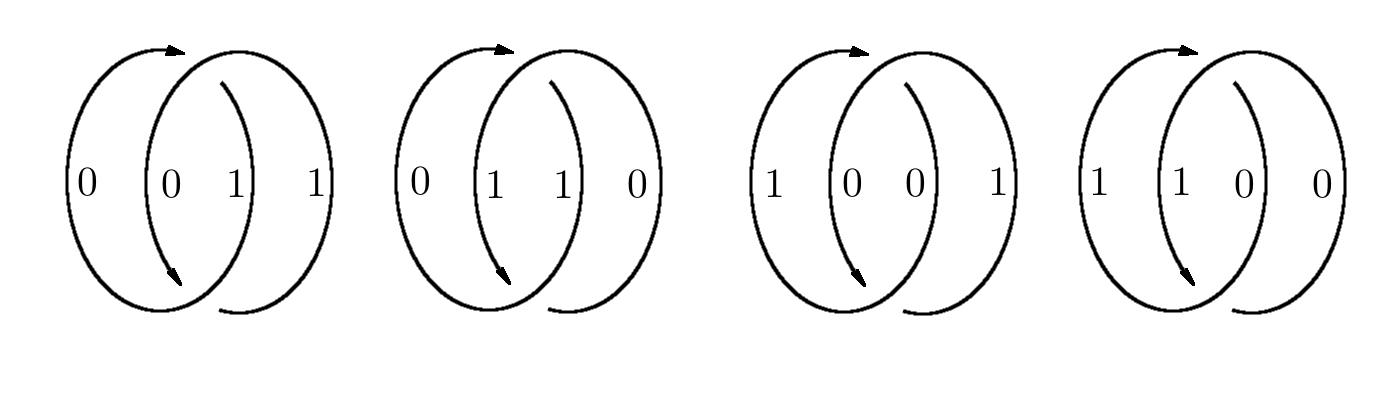}\]
The \textit{virtual Hopf link}, on the other hand, 
\[\includegraphics{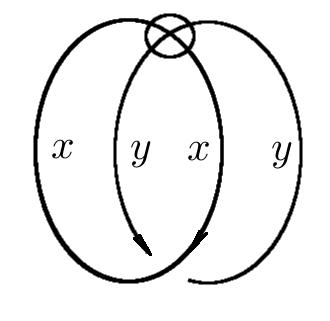} \quad 
\raisebox{0.5in}{$\begin{array}{rcl}
x & = & x+1 \\
y & = & y+1
\end{array}$}\]
has no colorings by $X$ since there are no solutions to the system
of coloring equations, so we have 
$\Phi_X^{\mathbb{Z}}(\mathrm{Hopf})=4\ne0=\Phi_X^{\mathbb{Z}}(\mathrm{vHopf})$
and the counting invariant distinguishes these two virtual links.
\end{example}

\section{\large\textbf{Biquandle Virtual Brackets}}\label{VB}

The invariant in Example \ref{ex:HL} gives a value of 
\[\Phi_X^{\mathbb{Z}}(L)=\left\{\begin{array}{ll} 2^c, & L\not\in\mathcal{O}  \\
0, & L\in\mathcal{O} \end{array}\right.\]
where $\mathcal{O}$ is the set of virtual links in which any component has
an odd number of (classical) crossing points and $c$ is the number of 
components of $L$.  Thus $\Phi_X^{\mathbb{Z}}$ with the choice of biquandle $X$ 
in example \ref{ex:HL} fails to distinguish any two classical knots; we would 
like to strengthen and \textit{enhance} this and other biquandle counting 
invariants to obtain stronger invariants. If $\beta$ is an invariant of 
$X$-colored 
diagrams, then the multiset of $\beta$-values over the set of $X$-colorings
of $L$ is a generally stronger invariant from which we can recover
$\Phi_X^{\mathbb{Z}}$ by taking the multiset's cardinality. In this section we 
introduce a new infinite family of such enhancements.

\begin{definition}
Let $X$ be a biquandle and $R$ a commutative ring with identity. Then a 
\textit{biquandle virtual bracket} consists of six maps 
$A,B,C,D,U,V:X\times X\to R$ and two distinguished elements
$\delta\in R$, $w\in R^{\times}$ satisfying the following conditions:
\[
\begin{array}{rcll}
w & = & \delta A_{x,x}+B_{x,x}+V_{x,x}, & (i.i)\\
w^{-1} & = &\delta C_{x,x}+D_{x,x}+U_{x,x}, & (i.ii),
\end{array}
\]
\[
\begin{array}{rcll}
1 & = & A_{x,y}C_{x,y}+V_{x,y}U_{x,y}, & (ii.i) \\
1 & = & B_{x,y}D_{x,y}+V_{x,y}U_{x,y}, & (ii.ii)\\
0 & = & A_{x,y}U_{x,y} +V_{x,y}C_{x,y}, & (ii.iii)\\
0 & = & B_{x,y}U_{x,y} +V_{x,y}D_{x,y}, & (ii.iv)\\
0 & = & \delta B_{x,y}D_{x,y}+A_{x,y}D_{x,y}+B_{x,y}C_{x,y}, & (ii.v)\\
0 & = & \delta A_{x,y}C_{x,y}+A_{x,y}D_{x,y}+B_{x,y}C_{x,y}, & (ii.vi),
\end{array}
\]
\[
\begin{array}{rcll}
A_{x,y}A_{x^y,z_y}A_{y,z}+V_{x,y}A_{x^y,z_y}V_{y,z} & = &A_{y_x,z_x}A_{x,z}A_{x^z,y^z}+V_{y_x,z_x}A_{x,z}V_{x^z,y^z} & (iii.i)\\
A_{x,y}A_{x^y,z_y}B_{y,z}+B_{x,y}A_{x^y,z_y}A_{y,z}  & & \\
+\delta B_{x,y}A_{x^y,z_y}B_{y,z}+ B_{x,y}A_{x^y,z_y}V_{y,z} & &\\
+B_{x,y}B_{x^y,z_y}B_{y,z}  +B_{x,y}V_{x^y,z_y}B_{y,z} & & \\
+V_{x,y}A_{x^y,z_y}B_{y,z} & = & A_{y_x,z_x}B_{x,z}A_{x^z,y^z} & (iii.ii)\\
A_{x,y}B_{x^y,z_y}A_{y,z} & = & A_{y_x,z_x}A_{x,z}B_{x^z,y^z}+ B_{y_x,z_x}A_{x,z}A_{x^z,y^z} \\
& & +\delta B_{y_x,z_x}A_{x,z}B_{x^z,y^z}+ B_{y_x,z_x}A_{x,z}V_{x^z,y^z}\\
& & + B_{y_x,z_x}B_{x,z}B_{x^z,y^z}+ B_{y_x,z_x}V_{x,z}B_{x^z,y^z}+ V_{y_x,z_x}A_{x,z}B_{x^z,y^z}
& (iii.iii)\\
A_{x,y}V_{x^y,z_y}A_{y,z} & = & A_{y_x,z_x}A_{x,z}V_{x^z,y^z}+V_{y_x,z_x}A_{x,z}A_{x^z,y^z} & (iii.iv)\\
A_{x,y}A_{x^y,z_y}V_{y,z}+V_{x,y}A_{x^y,z_y}A_{y,z} & = & A_{y_x,z_x}V_{x,z}A_{x^z,y^z} & (iii.v)\\
B_{x,y}B_{x^y,z_y}A_{y,z}+B_{x,y}V_{x^y,z_y}V_{y,z} & = & A_{y_x,z_x}B_{x,z}B_{x^z,y^z}+ V_{y_x,z_x}V_{x,z}B_{x^z,y^z} & (iii.vi)\\
A_{x,y}B_{x^y,z_y}B_{y,z}+V_{x,y}V_{x^y,z_y}B_{y,z} & = & B_{y_x,z_x}B_{x,z}A_{x^z,y^z}+ B_{y_x,z_x}V_{x,z}V_{x^z,y^z} & (iii.vii)\\
B_{x,y}B_{x^y,z_y}V_{y,z}+B_{x,y}V_{x^y,z_y}A_{y,z} & = & A_{y_x,z_x}B_{x,z}V_{x^z,y^z} & (iii.viii)\\
A_{x,y}B_{x^y,z_y}V_{y,z} & = & B_{y_x,z_x}B_{x,z}V_{x^z,y^z}+B_{y_x,z_x}V_{x,z}A_{x^z,y^z}& (iii.ix)\\
V_{x,y}B_{x^y,z_y}A_{y,z} & = & A_{y_x,z_x}V_{x,z}B_{x^z,y^z}+V_{y_x,z_x}B_{x,z}B_{x^z,y^z} & (iii.x)\\
A_{x,y}V_{x^y,z_y}B_{y,z}+V_{x,y}B_{x^y,z_y}B_{y,z} & = & V_{y_x,z_x}B_{x,z}A_{x^z,y^z} & (iii.xi)\\
V_{x,y}V_{x^y,z_y}A_{y,z} & = & A_{y_x,z_x}V_{x,z}V_{x^z,y^z} & (iii.xii)\\
A_{x,y}V_{x^y,z_y}V_{y,z} & = & V_{y_x,z_x}V_{x,z}A_{x^z,y^z} & (iii.xiii)\\
V_{x,y}B_{x^y,z_y}V_{y,z} & = & V_{y_x,z_x}B_{x,z}V_{x^z,y^z} & (iii.xiv)\\
V_{x,y}V_{x^y,z_y}V_{y,z} & = & V_{y_x,z_x}V_{x,z}V_{x^z,y^z} & (iii.xv)
\end{array}\]
Note that in the case that $A_{x,y},B_{x,y},C_{x,y}$ and $D_{x,y}$ are all
invertible for a pair $(x,y)$, conditions (ii.iii) and (ii.iv) reduce 
to 
\[\begin{array}{rcll}
\delta & = & -A_{x,y}B_{x,y}^{-1}-C_{x,y}D_{x,y}^{-1}, & (ii.iii')\\
\delta & = & -A_{x,y}^{-1}B_{x,y}-C_{x,y}^{-1}D_{x,y} & (ii.iv').
\end{array}\]
\end{definition}


The biquandle virtual bracket axioms are chosen so that the skein relations
\[\includegraphics{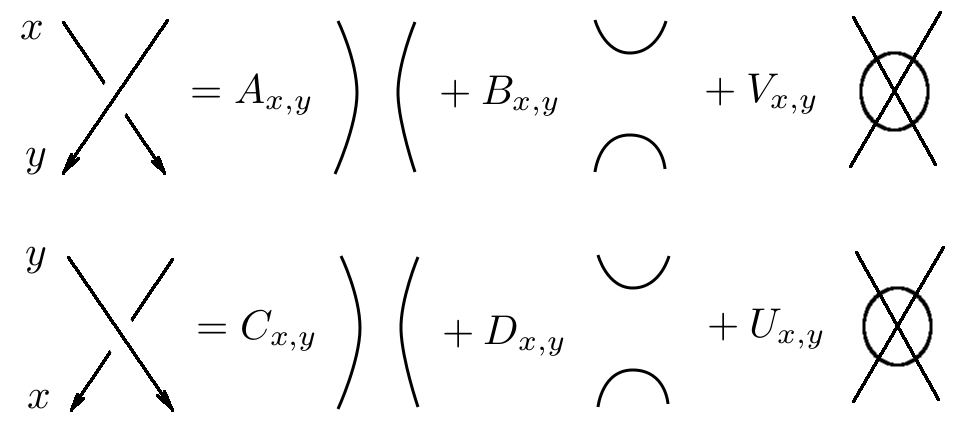}\]
with writhe correction factor $w$ and smoothed component value $\delta$
\[\includegraphics{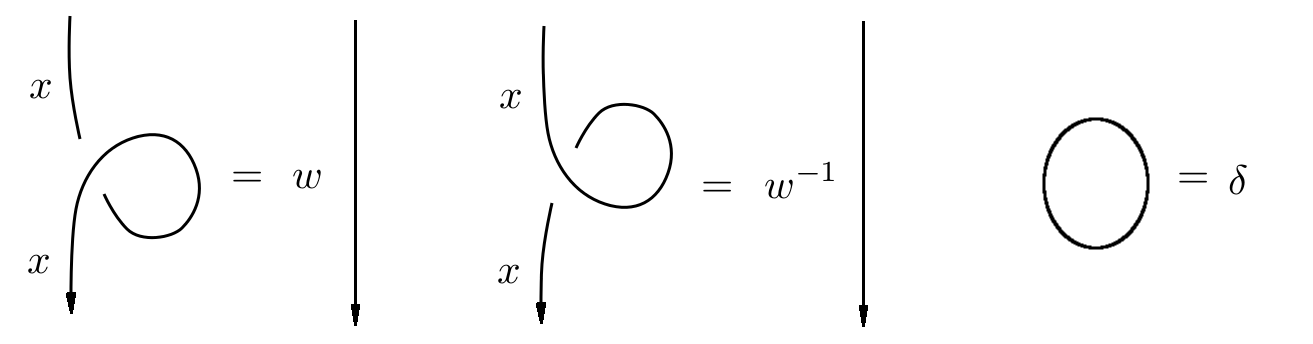}\]
define an invariant of $X$-colored Reidemeister moves.

To accomplish this, we start with the idea that for a given $X$-coloring $f$ 
of an oriented link diagram $L$ with $c$ classical crossings, there is a set of 
$3^c$ completely smoothed
diagrams known as \textit{states} in which each classical crossing is 
smoothed either with the orientation, against the orientation, or made virtual.
Each such state $S$ has an associated value $\beta_S\in R$ defined as the 
product of the coefficients associated to the choices of smoothings times 
$\delta^{k}$ where $k$ is the number of closed curve components of the state
(which may have virtual crossings) times $w^{n-p}$ where $n$ is the number of 
negative classical crossings and $p$ is the number of positive classical 
crossings of $L$. We then sum these $\beta_S$ values over the complete set
of states obtained from $f$ to obtain the \textit{state sum value} $\beta_f$
for the given coloring of the link, and the multiset of $\beta$ values over the 
set of all $X$-colorings $f$ of $L$ is then an invariant of oriented links.

We have the following critical observation:

\begin{observation}\label{ob1}
A move on a tangle which does not change connectivity on the boundary
(i.e., boundary points connected or not connected in the pre-move tangle 
are still connected or not connected in the post-move tangle) and for which
the coefficient weight sums are equal before and after the move does not 
change the overall contribution to the state sum.
\end{observation}

Using this observation, we obtain the biquandle virtual bracket axioms 
from the oriented virtual Reidemeister moves.

The Reidemeister I moves require that
$w =  \delta A_{x,x}+B_{x,x}+V_{x,x}$ and
$w^{-1} = \delta C_{x,x}+D_{x,x}+U_{x,x}$:
\[\includegraphics{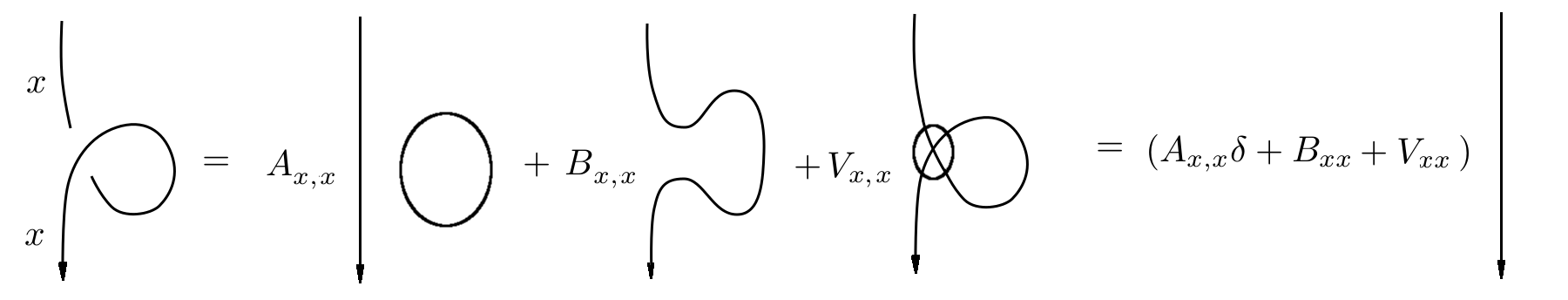}\]
\[\includegraphics{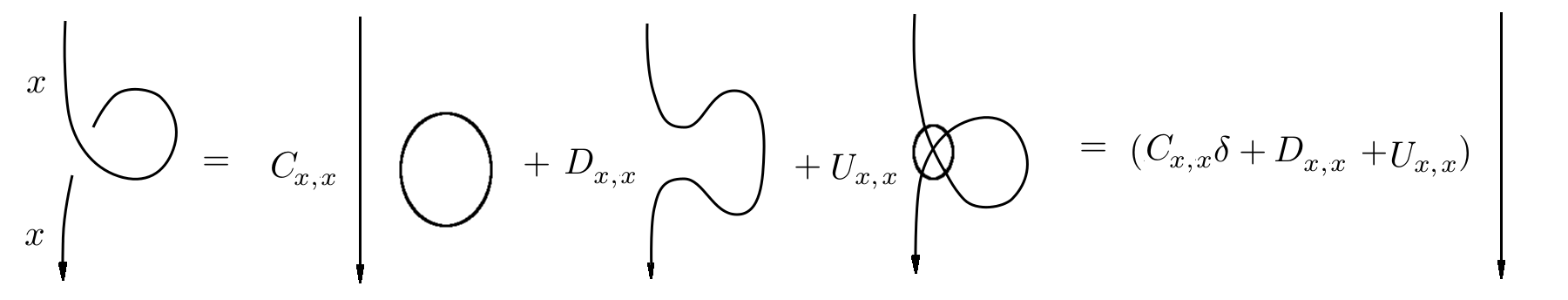}\]

Comparing the skein expansions before and after the direct and reverse
oriented Reidemeister II moves yields the conditions
\[\begin{array}{rcl}
1 & = & A_{x,y}C_{x,y}+V_{x,y}U_{x,y} \\
1 & = & B_{x,y}D_{x,y}+V_{x,y}U_{x,y} \\
0 & = & A_{x,y}U_{x,y} +V_{x,y}C_{x,y} \\
0 & = & B_{x,y}U_{x,y} +V_{x,y}D_{x,y} \\
0 & = & \delta B_{x,y}D_{x,y}+A_{x,y}D_{x,y}+B_{x,y}C_{x,y}\\
0 & = & \delta A_{x,y}C_{x,y}+A_{x,y}D_{x,y}+B_{x,y}C_{x,y}\\
\end{array}\]

\[\includegraphics{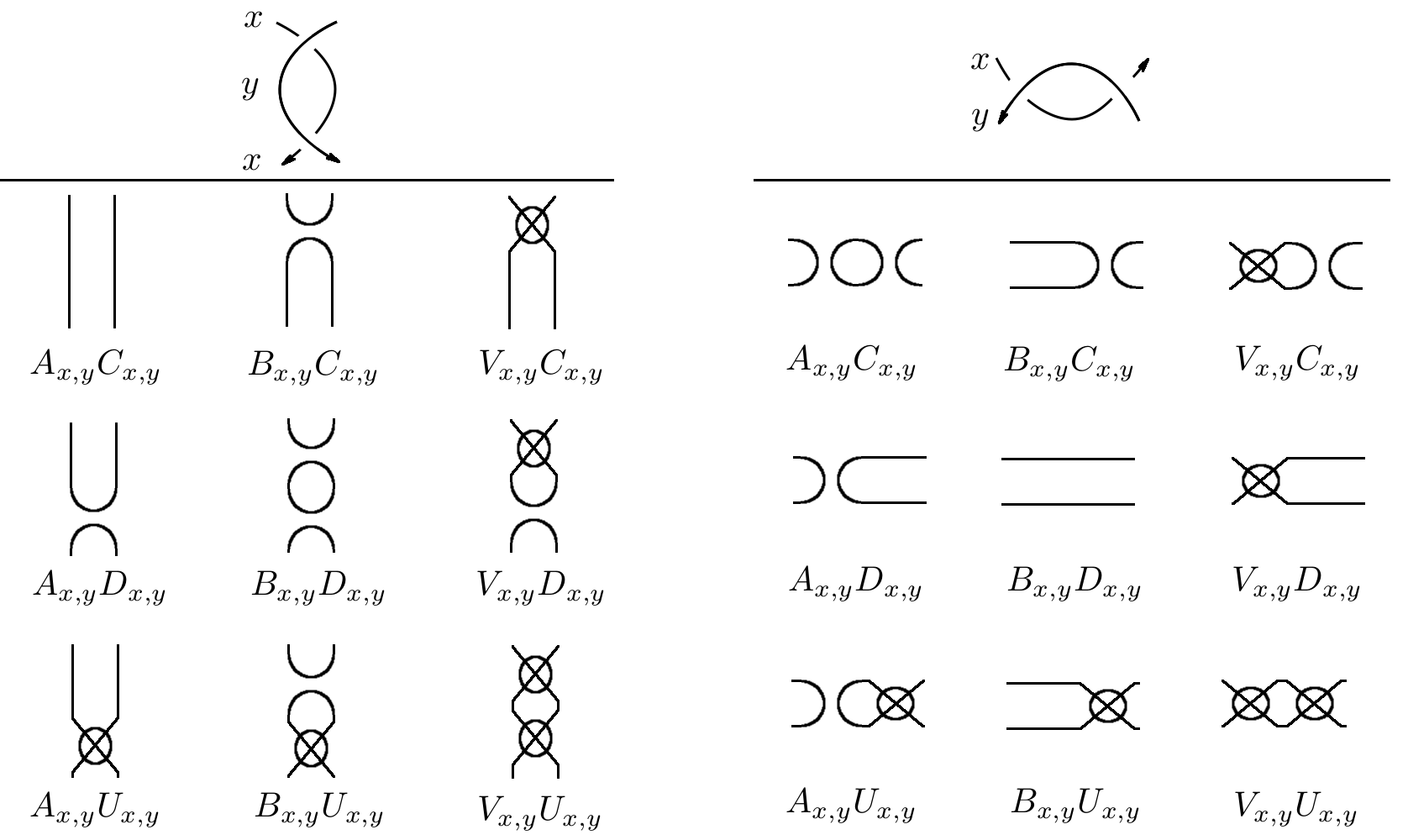}\]

The two sides of the third Reidemeister move, with each smoothed state
labeled with its coefficient product,  
\[\includegraphics{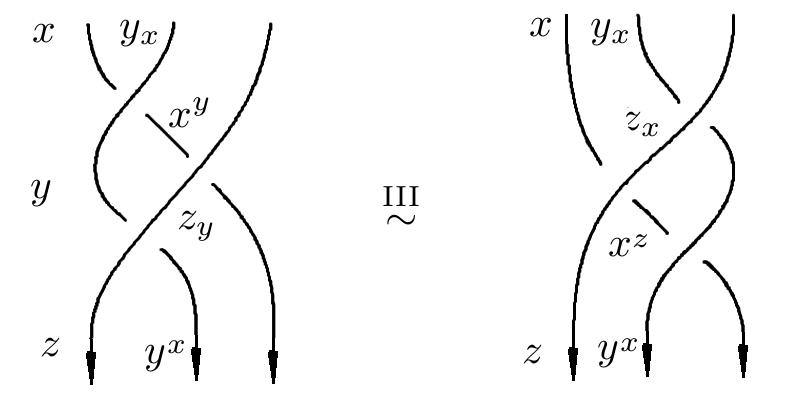}\] 
are left-hand side
\[\includegraphics{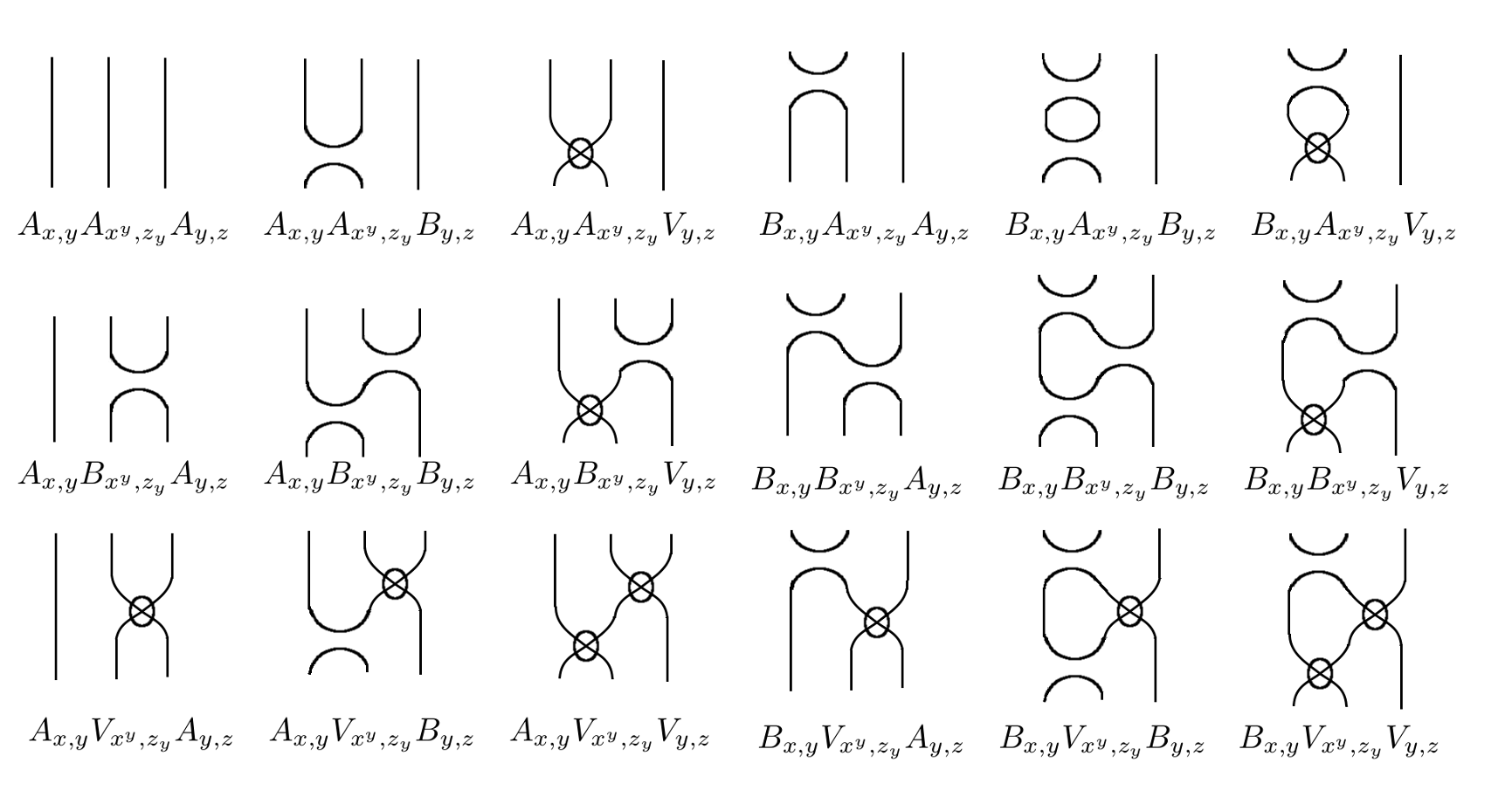}\]
\[\includegraphics{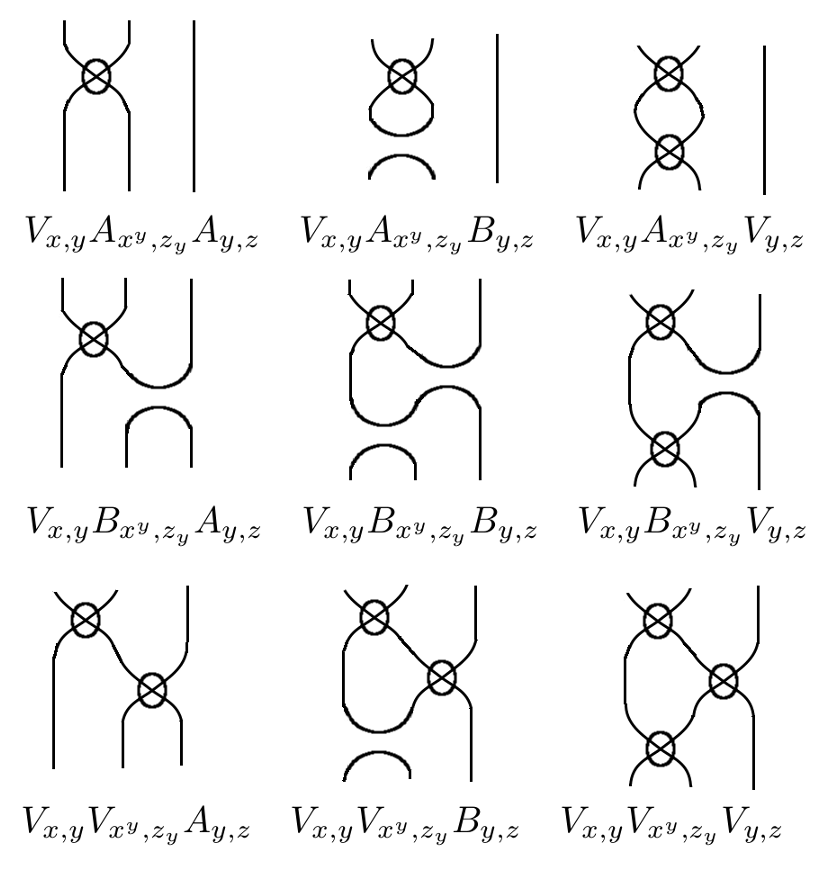}\]
and right-hand side
\[\includegraphics{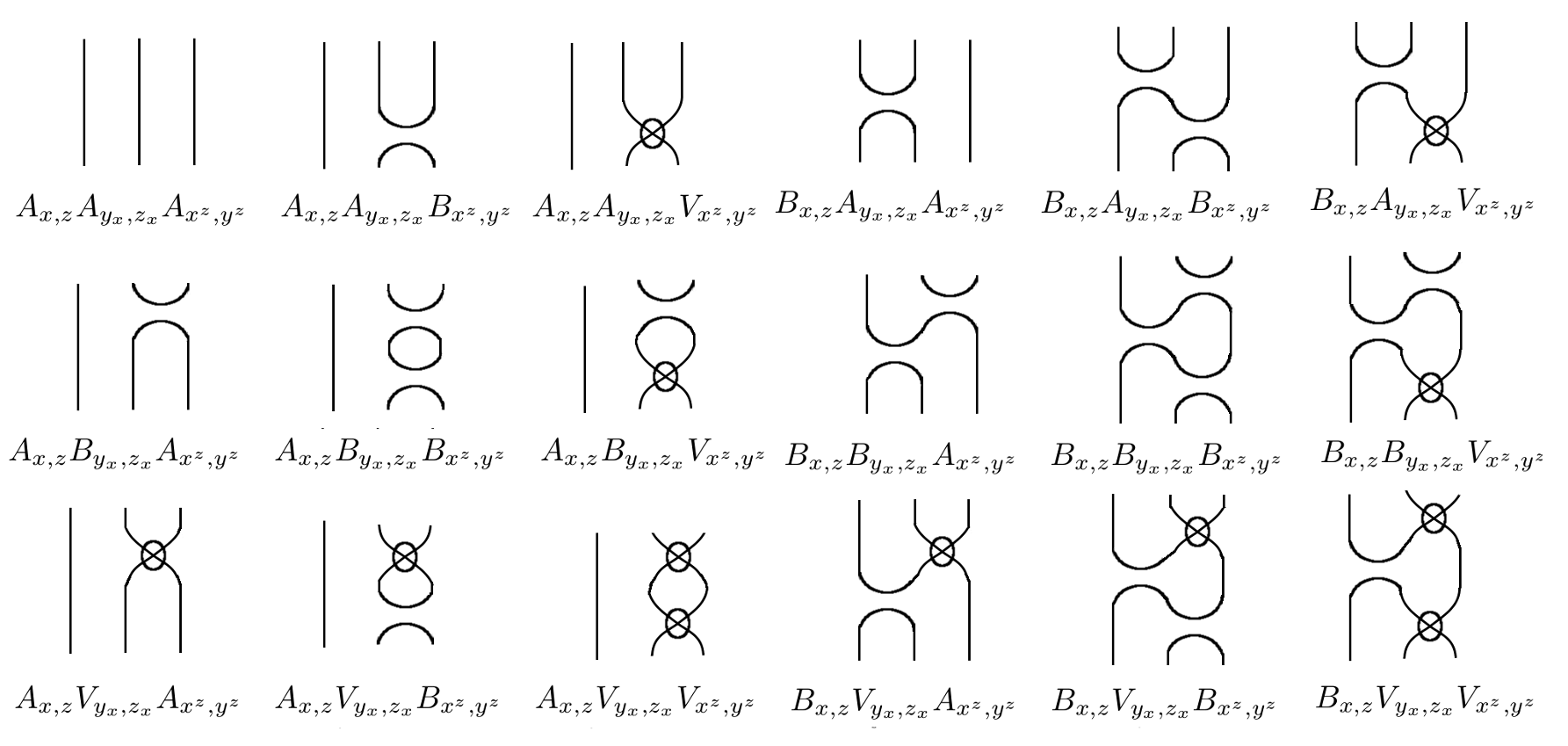}\]
\[\includegraphics{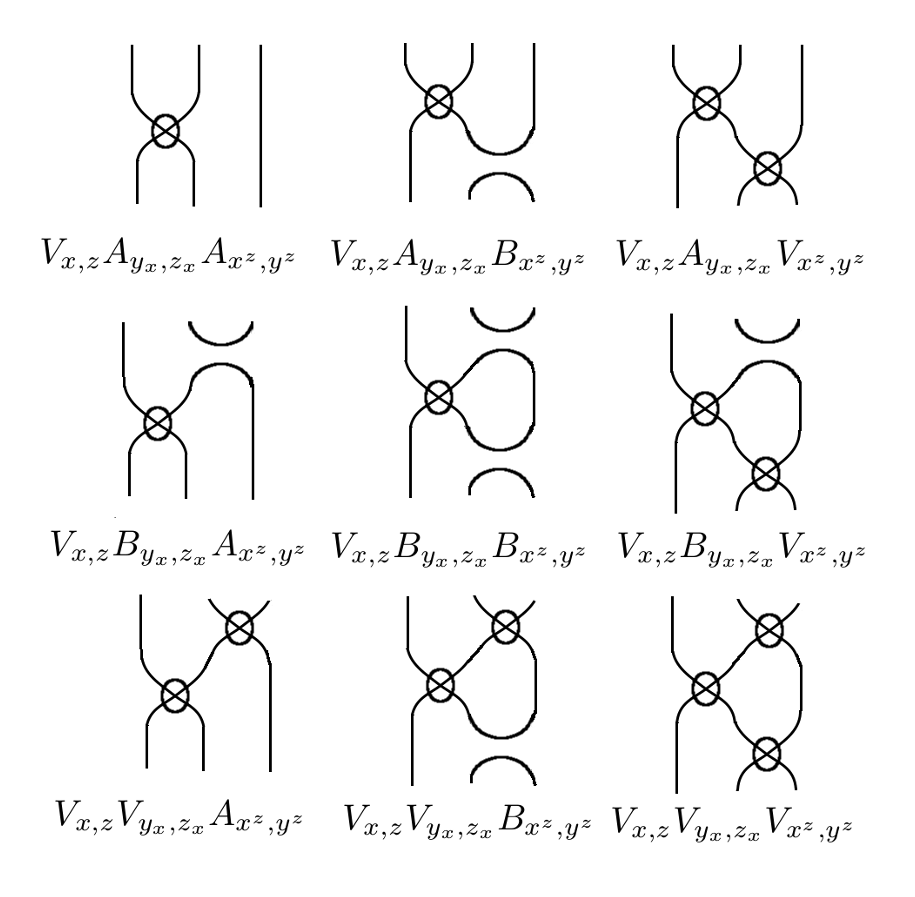}.\]
Comparing coefficient sums of equivalent tangles using observation \ref{ob1}
yields the fifteen axioms $(iii.i)$ through $(iii.xv)$.

The pure virtual moves do not affect the state sum value in any way, and the 
mixed virtual move yields the same sum on both sides up to virtual moves, and
thus does not change the state sum:
\[\includegraphics{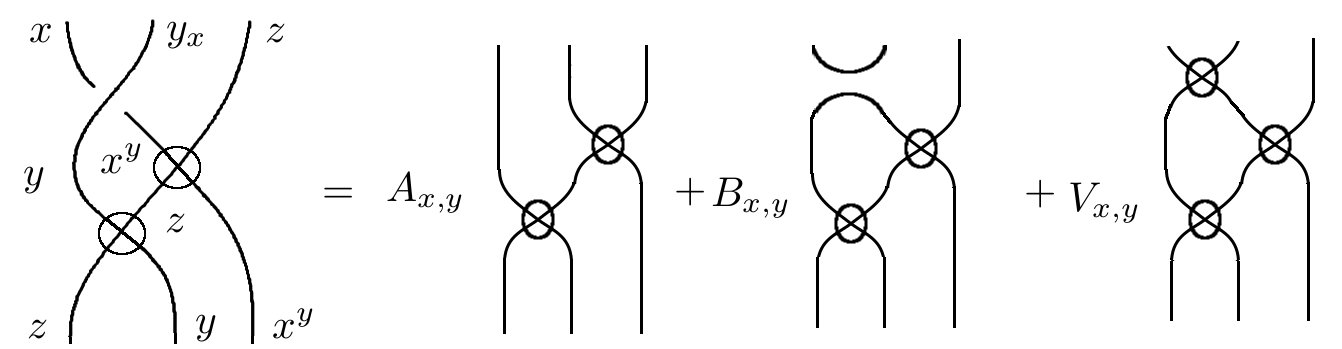}\]
\[\includegraphics{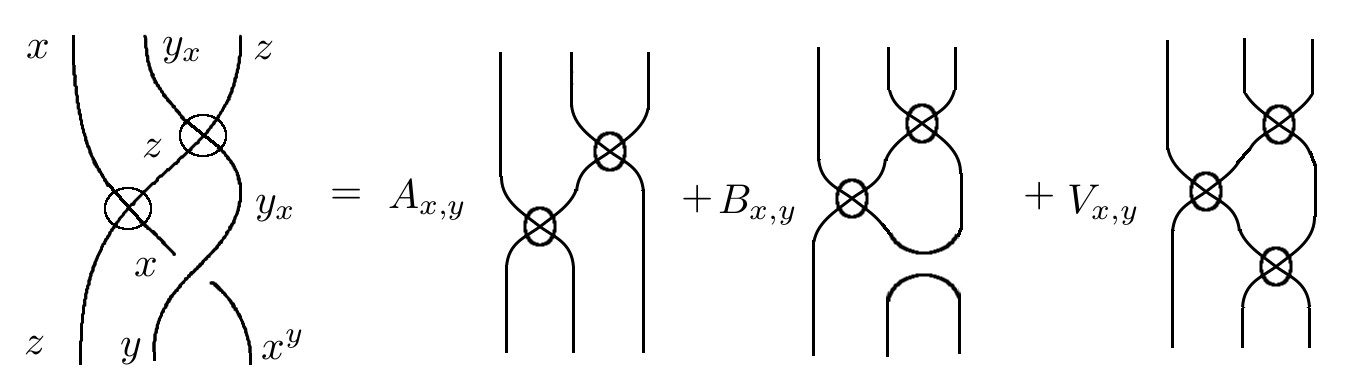}\]
Hence, we have
\begin{theorem}
Let $X$ be a biquandle and $\beta$ a biquandle virtual bracket. Then for any 
$X$-coloring $f$ of an oriented (classical or) virtual link $L$, the state sum
\[\beta_f=\sum_{S\ \mathrm{State}} \beta_S\]
is an invariant of $X$-colored virtual Reidemeister moves.
\end{theorem}

\begin{example}
Treating an oriented virtual link $L$ as colored by its fundamental biquandle
$\mathcal{B}(L)$ with $f$ assigning to each semiarc its generator, the state sum
$\beta_f$ is the \textit{fundamental biquandle virtual bracket} of $L$. While
it is difficult to compare values of the fundamental virtual bracket directly
(since doing so requires comparing fundamental biquandles of links directly), we
find it useful to use the fundamental virtual bracket expression as a step in 
computing state sums for virtual brackets over finite biquandles by 
substituting image values of generators under the coloring maps 
$f\in\mathrm{Hom}(\mathcal{B}(L),X)$; see Section \ref{AE}.
\end{example}

\begin{corollary}\label{cor:inv}
For any biquandle $X$, biquandle virtual bracket $\beta$ and oriented virtual 
link $L$, the multiset of state sum values
\[\Phi_X^{\beta,M}(L)=\{\beta_f\ |\ f\in\mathrm{Hom}(\mathcal{B}(L),X)\}\]
is an invariant of classical and virtual links.
\end{corollary}

\begin{definition}
We call the multiset in Corollary \ref{cor:inv} the \textit{biquandle virtual 
bracket multiset} invariant of $L$ with respect to the biquandle $X$ and
the biquandle virtual bracket $\beta$. 
\end{definition}

If $R$ is a number ring, it is common practice (see \cite{CJKLS,EN} etc.) to
write the multiset in a ``polynomial'' form vaguely analogous to a generating
function by making the multiplicities coefficients and the elements exponents,
e.g. encoding $\{-2,1,1,1,3,3\}$ as $u^{-2}+3u+2u^3$. When written this way, 
we will refer to the invariant as the \textit{biquandle virtual bracket 
polynomial} of $L$, denoted $\Phi_X^{\beta}(L)$.

\begin{example}
Let $X$ be a biquandle and $\beta$ a biquandle virtual bracket over a 
commutative ring with identity $R$. If $V_{x,y}=U_{x,y}=0$ for all $x,y\in X$,
then $\beta$ is a biquandle bracket as defined in \cite{NOR}. We will
refer to such virtual brackets as \textit{biquandle classical brackets}
or just \textit{biquandle brackets}.
\end{example}

\begin{example}
Let $X$ be a biquandle and $\beta$ a biquandle virtual bracket over a 
commutative ring with identity $R$. If $A_{x,y}=B_{x,y}=C_{x,y}=D_{x,y}=0$ for
all $x,y\in X$ and $w=\delta=1$, then the virtual bracket axioms reduce to
\[U_{x,y}=V_{x,y}^{-1}\quad \mathrm{and}\quad 
V_{x,y}V_{x^y,z_y}V_{y,z} = V_{y_x,z_x}V_{x,z}V_{x^z,y^z}\]
for all $x,y,z\in X$, making $V$ a biquandle 2-cocycle and $\Phi_X^{B}$ 
exactly the biquandle 2-cocycle invariant with Boltzmann weights written 
multiplicatively (see \cite{EN,CEGN,CES}).
\end{example}

\section{Applications and Examples}\label{AE}
In this section, we provide examples which show that our invariants, 
including the special case of biquandle classical brackets from \cite{NOR}, 
can be used to determine chirality of a virtual knot or link, or to determine 
non-invertibility of a virtual knot or link.

\begin{example}\label{ex:vtrefoil}
Let $X=\mathbb Z_2$ be the biquandle with operations 
$\underline{\triangleright}, \overline{\triangleright}: X \times X \to X$ 
defined by
\[
x\,\underline{\triangleright}\,y= x+1 \mbox{\ \ \  and \ \ \ } x\,\overline{\triangleright}\,y= x+1.
\]
Let $\delta=3 \in \mathbb Z_5$.
Define $A_{xy}, B_{xy}, V_{xy}, C_{xy}, D_{xy}, U_{xy}: X \times X \to \mathbb Z_5$ by
\[
[A_{xy}| B_{xy}| V_{xy}| C_{xy}| D_{xy}| U_{xy}]=
\left[
\begin{array}{cc|cc|cc|cc|cc|cc}
1&1&1&1&0&2&1&3&1&3&0&4\\
1&1&1&1&3&0&3&1&3&1&1&0
\end{array}
\right].
\] 
Our Python and Mathematica computations agree that this data defines a 
biquandle virtual bracket structure over $X$. 

The right-handed virtual trefoil  $L_1^{\rm vtrefoil}$ below has nine states with 
coefficients as listed.
\[  \includegraphics[clip,width=7.0cm]{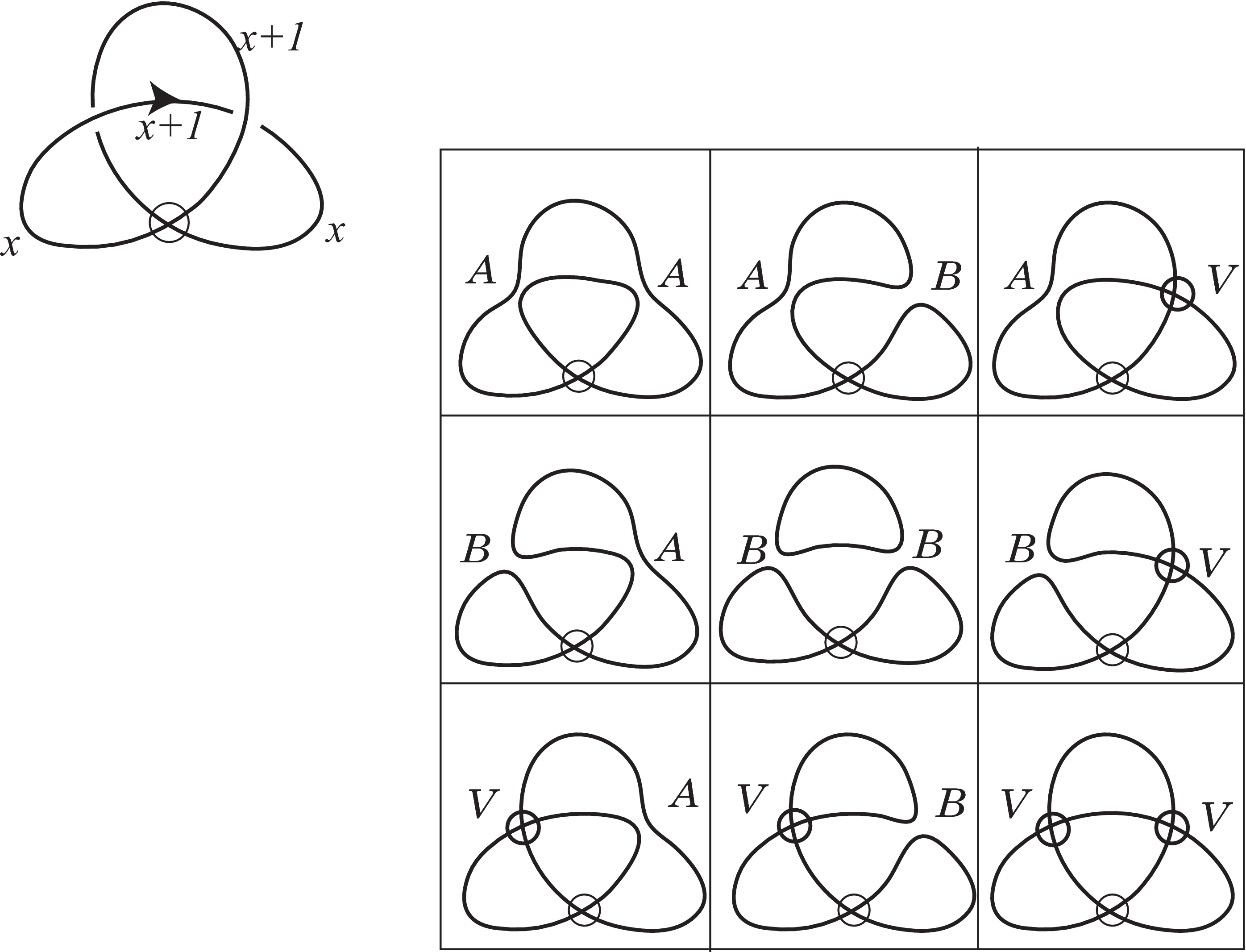}\]
Then the fundamental biquandle virtual bracket value with respect to the biquandle $X$ is 
\[
\begin{array}{l}
\phi_1(x)
 = (A_{x\, x+1}V_{x+1\, x }+B_{x\, x+1}B_{x+1\, x  }+V_{x\, x+1}A_{x+1\, x  })\delta^2\\ 
\quad\quad\quad\quad
+(A_{x\, x+1}A_{x+1\, x }+A_{x\, x+1}B_{x+1\, x  }+B_{x\, x+1}A_{x+1\, x  }+B_{x\, x+1}V_{x+1\, x  }+V_{x\, x+1}B_{x+1\, x  }+V_{x\, x+1}V_{x+1\, x })\delta.
\end{array}
\]
Hence we have 
\[
\begin{array}{c|c}
x&\phi_1 (x)\\ \hline
0&(1\cdot3+1\cdot1+2\cdot1)\cdot3^2 +(1\cdot1+1\cdot1+1\cdot1+1\cdot3+2\cdot1+2\cdot3)\cdot3=1\\
1&(1\cdot2+1\cdot1+3\cdot1)\cdot3^2 +(1\cdot1+1\cdot1+1\cdot1+1\cdot2+3\cdot1+3\cdot2)\cdot3=1
\end{array}
\]
Therefore the right-handed   virtual trefoil  $L_1^{\rm vtrefoil}$ has biquandle virtual bracket invariant
\[
\Phi(L_1^{\rm vtrefoil})=2 u.
\]

The left-handed virtual trefoil $L_2^{\rm vtrefoil}$ below has nine states with coefficients as listed.
\[  \includegraphics[clip,width=7.0cm]{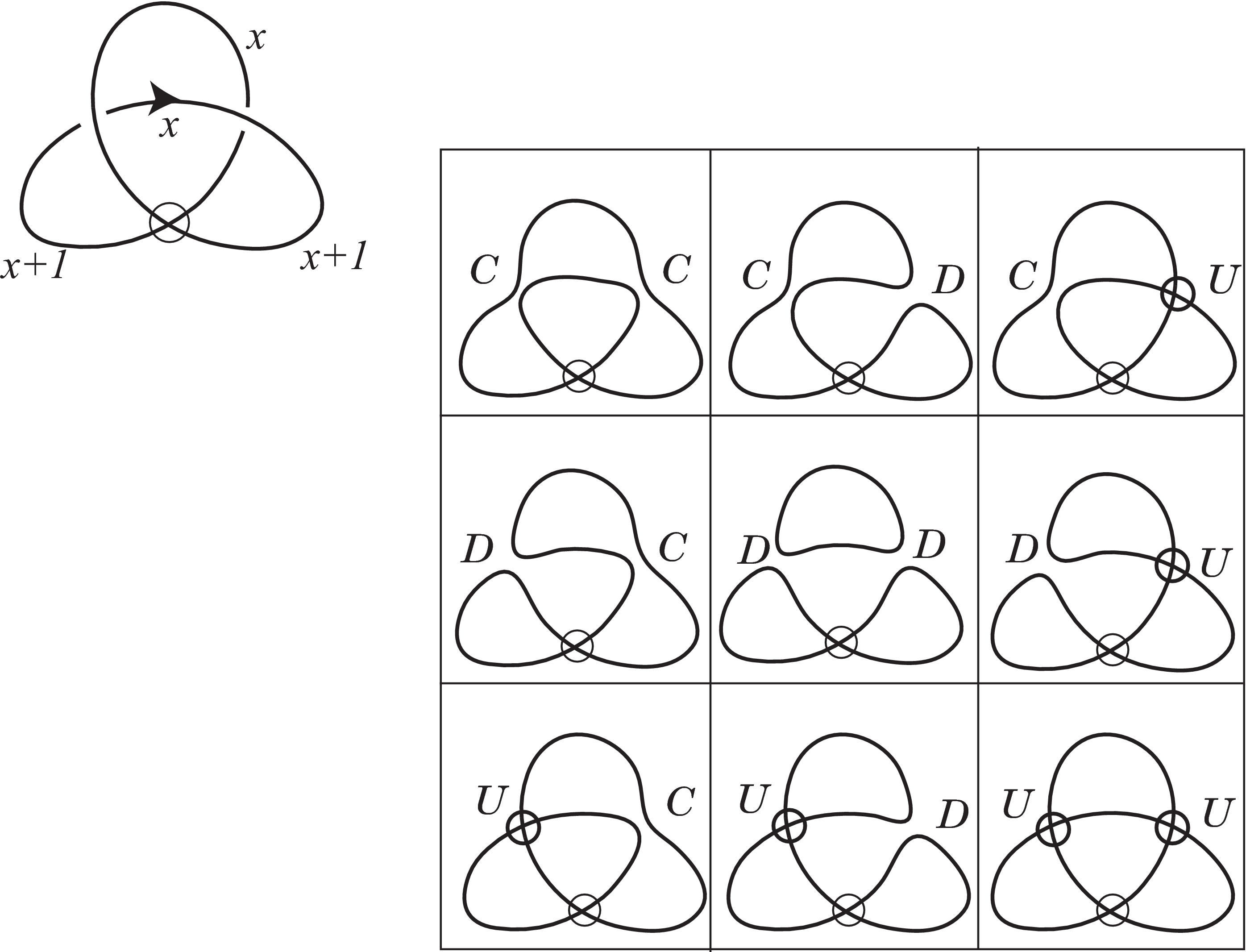}\]
Then the fundamental biquandle virtual bracket value with respect to $X$ is 
\[
\begin{array}{l}
\phi_2(x)=(C_{x\, x+1}U_{x+1\, x}+D_{x\, x+1}D_{x+1\, x}+U_{x\, x+1}C_{x+1\, x})\delta^2\\
\quad\quad\quad\quad
+(C_{x\, x+1}C_{x+1\, x}+C_{x\, x+1}D_{x+1\, x}+D_{x\, x+1}C_{x+1\, x}+D_{x\, x+1}U_{x+1\, x}+U_{x\, x+1}D_{x+1\, x}+U_{x\, x+1}U_{x+1\, x})\delta.
\end{array}
\]
Hence we have 
\[
\begin{array}{c|c}
x&\phi_2(x)\\ \hline
0&(3\cdot1+3\cdot3+4\cdot3)\cdot3^2 +(3\cdot3+3\cdot3+3\cdot3+3\cdot1+4\cdot3+4\cdot1)\cdot3=4\\
1&(3\cdot4+3\cdot3+1\cdot3)\cdot3^2 +(3\cdot3+3\cdot3+3\cdot3+3\cdot4+1\cdot3+1\cdot4)\cdot3=4
\end{array}
\]
Therefore the left-handed virtual trefoil  $L_2^{\rm vtrefoil}$ has biquandle 
virtual bracket invariant
\[
\Phi(L_2^{\rm vtrefoil})=2 u^{4}.
\]

Therefore  the biquandle virtual bracket invariants of the right- and 
left-handed virtual trefoils (with any orientation) have the different 
values, and thus, they are not equivalent. Note that it is known that the 
right-handed (or left-handed)  virtual trefoil  is equivalent to its inverse.
\end{example}


\begin{example}
Let $X=\mathbb Z_3$ be the biquandle with operations 
$ \underline{\triangleright}, \overline{\triangleright}: X \times X \to X$ 
defined by
\[
x\,\underline{\triangleright}\,y= x+2 \mbox{\ \ \  and \ \ \ } x\,\overline{\triangleright}\,y= x+2.
\]
Let $\delta=2 \in \mathbb Z_3$.
Define $A_{xy}, B_{xy}, V_{xy}, C_{xy}, D_{xy}, U_{xy}: X \times X \to \mathbb Z_3$ by
\[\begin{array}{l}
[A_{xy}| B_{xy}| V_{xy}| C_{xy}| D_{xy}| U_{xy}]\\[10pt]
=
\left[
\begin{array}{ccc|ccc|ccc|ccc|ccc|ccc}
1&1&2&2&2&1&0&0&0&1&1&2&2&2&1&0&0&0\\
2&1&2&1&2&1&0&0&0&2&1&2&1&2&1&0&0&0\\
1&1&1&2&2&2&0&0&0&1&1&1&2&2&2&0&0&0
\end{array}
\right].
\end{array}
\] 
Our Python and Mathematica computations agree that this defines a biquandle 
virtual bracket structure. 

Let $K_1$ be the oriented virtual-knot depicted below.  
\[  \includegraphics[clip,width=5.0cm]{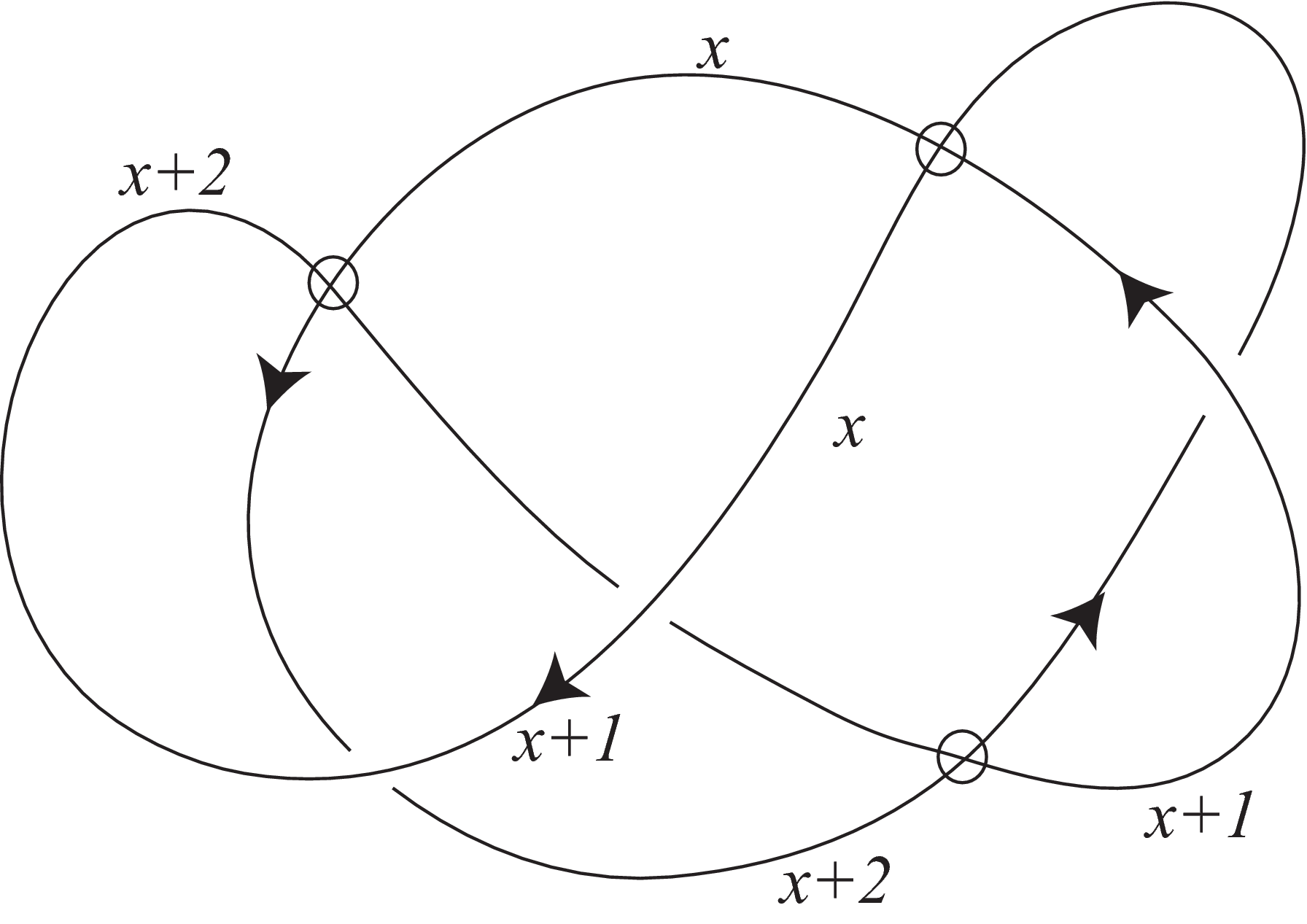}\]
$K_1$ has twenty-seven states 
and 
its fundamental biquandle virtual bracket value with respect to the biquandle $X$ is 
\[
\begin{array}{ll}
\phi_1 (x)&=V_{x\, x+2}A_{x+2\, x+1}C_{x\, x+1}\delta^3\\
&+(A_{x\, x+2}A_{x+2\, x+1}C_{x\, x+1}+A_{x\, x+2}B_{x+2\, x+1}D_{x\, x+1}\\
&+A_{x\, x+2}V_{x+2\, x+1}U_{x\, x+1}+B_{x\, x+2}A_{x+2\, x+1}C_{x\, x+1}\\
&+B_{x\, x+2}B_{x+2\, x+1}U_{x\, x+1}+B_{x\, x+2}V_{x+2\, x+1}D_{x\, x+1}\\
&+V_{x\, x+2}A_{x+2\, x+1}D_{x\, x+1}+V_{x\, x+2}A_{x+2\, x+1}U_{x\, x+1}\\
&+V_{x\, x+2}B_{x+2\, x+1}C_{x\, x+1}+V_{x\, x+2}V_{x+2\, x+1}C_{x\, x+1})\delta^2\\
&+( A_{x\, x+2}A_{x+2\, x+1}D_{x\, x+1}+A_{x\, x+2}A_{x+2\, x+1}U_{x\, x+1}\\
&+A_{x\, x+2}B_{x+2\, x+1}C_{x\, x+1}+A_{x\, x+2}B_{x+2\, x+1}U_{x\, x+1}\\
&+A_{x\, x+2}V_{x+2\, x+1}C_{x\, x+1}+A_{x\, x+2}V_{x+2\, x+1}D_{x\, x+1}\\
&+B_{x\, x+2}A_{x+2\, x+1}D_{x\, x+1}+B_{x\, x+2}A_{x+2\, x+1}U_{x\, x+1}\\
&+B_{x\, x+2}B_{x+2\, x+1}C_{x\, x+1}+B_{x\, x+2}B_{x+2\, x+1}D_{x\, x+1}\\
&+B_{x\, x+2}V_{x+2\, x+1}C_{x\, x+1}+B_{x\, x+2}V_{x+2\, x+1}U_{x\, x+1}\\
&+V_{x\, x+2}B_{x+2\, x+1}D_{x\, x+1}+V_{x\, x+2}B_{x+2\, x+1}U_{x\, x+1}\\
&+V_{x\, x+2}V_{x+2\, x+1}D_{x\, x+1}+V_{x\, x+2}V_{x+2\, x+1}U_{x\, x+1})\delta,
\end{array}
\]
and we can check that $\phi_1(0)=\phi_1(1)=\phi_1(2)=1$.
Therefore $K_1$ has biquandle virtual bracket invariant
\[
\Phi(K_1)=3 u.
\]

Let $K_2$ be the oriented virtual-knot depicted below, that is,  
it is the same virtual-knot $K_1$ with its orientation reversed.
\[  \includegraphics[clip,width=5.0cm]{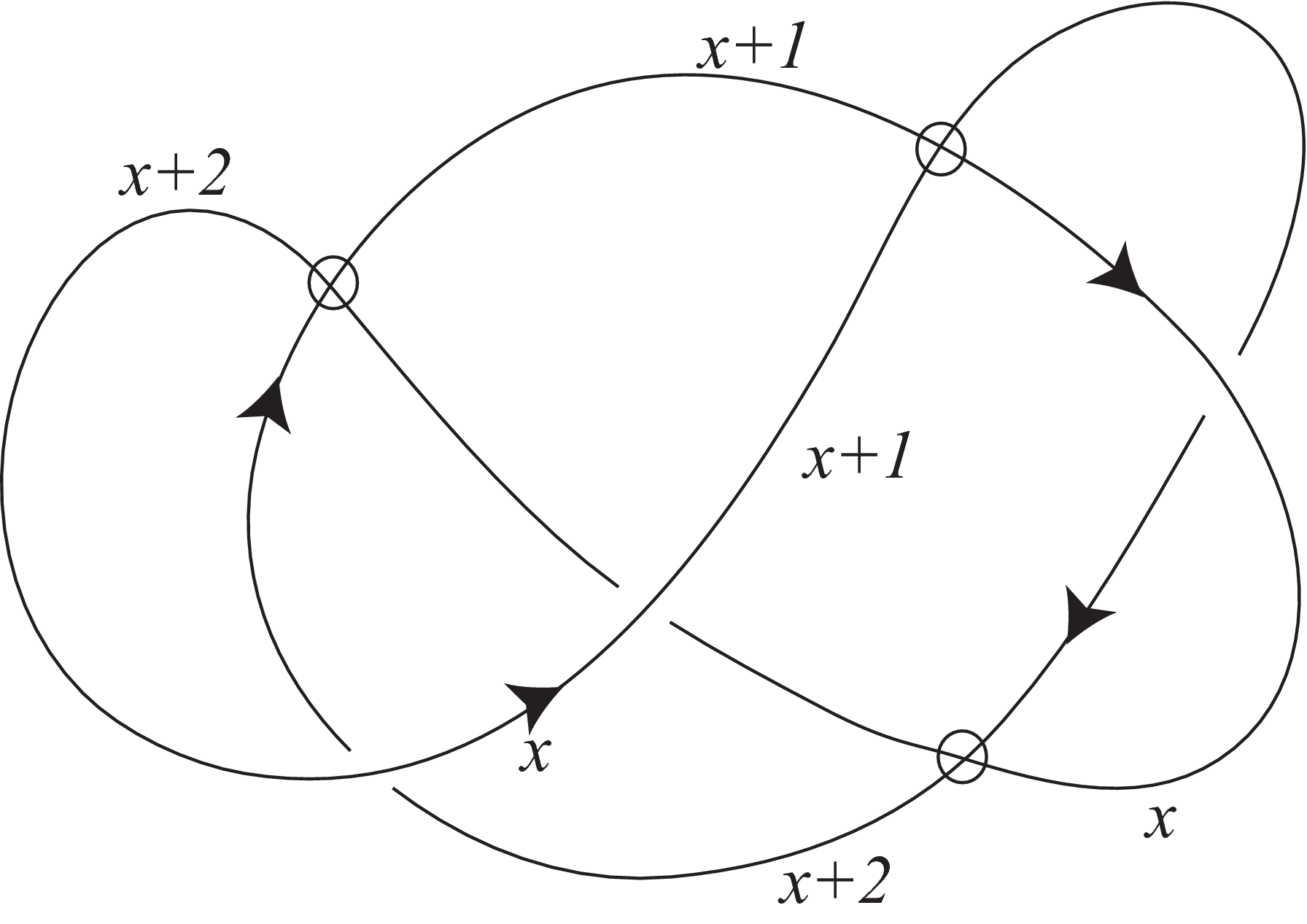}\]
$K_2$ has twenty-seven states 
and 
its fundamental biquandle virtual bracket value with respect to the biquandle $X$ is 
\[
\begin{array}{ll}
\phi_2 (x)&=V_{x+2\, x}A_{x\, x+1}C_{x+2\, x+1}\delta^3\\
&+(A_{x+2\, x}A_{x\, x+1}C_{x+2\, x+1}+A_{x+2\, x}B_{x\, x+1}D_{x+2\, x+1}\\
&+A_{x+2\, x}V_{x\, x+1}U_{x+2\, x+1}+B_{x+2\, x}A_{x\, x+1}C_{x+2\, x+1}\\
&+B_{x+2\, x}B_{x\, x+1}U_{x+2\, x+1}+B_{x+2\, x}V_{x\, x+1}D_{x+2\, x+1}\\
&+V_{x+2\, x}A_{x\, x+1}D_{x+2\, x+1}+V_{x+2\, x}A_{x\, x+1}U_{x+2\, x+1}\\
&+V_{x+2\, x}B_{x\, x+1}C_{x+2\, x+1}+V_{x+2\, x}V_{x\, x+1}C_{x+2\, x+1})\delta^2\\
&+( A_{x+2\, x}A_{x\, x+1}D_{x+2\, x+1}+A_{x+2\, x}A_{x\, x+1}U_{x+2\, x+1}\\
&+A_{x+2\, x}B_{x\, x+1}C_{x+2\, x+1}+A_{x+2\, x}B_{x\, x+1}U_{x+2\, x+1}\\
&+A_{x+2\, x}V_{x\, x+1}C_{x+2\, x+1}+A_{x+2\, x}V_{x\, x+1}D_{x+2\, x+1}\\
&+B_{x+2\, x}A_{x\, x+1}D_{x+2\, x+1}+B_{x+2\, x}A_{x\, x+1}U_{x+2\, x+1}\\
&+B_{x+2\, x}B_{x\, x+1}C_{x+2\, x+1}+B_{x+2\, x}B_{x\, x+1}D_{x+2\, x+1}\\
&+B_{x+2\, x}V_{x\, x+1}C_{x+2\, x+1}+B_{x+2\, x}V_{x\, x+1}U_{x+2\, x+1}\\
&+V_{x+2\, x}B_{x\, x+1}D_{x+2\, x+1}+V_{x+2\, x}B_{x\, x+1}U_{x+2\, x+1}\\
&+V_{x+2\, x}V_{x\, x+1}D_{x+2\, x+1}+V_{x+2\, x}V_{x\, x+1}U_{x+2\, x+1})\delta,
\end{array}
\]
and we can check that $\phi_1(0)=\phi_1(1)=\phi_1(2)=2$.
Therefore $K_2$ has biquandle bracket invariant
\[
\Phi(K_2)=3 u^2.
\]

Therefore  the biquandle virtual bracket invariants of $K_1$ and $K_2=-K_1$ have the different values, and thus, they are not equivalent. In particular,
$\Phi_X^{B}$ can detect invertibility. 
\end{example}



\begin{example} \label{ex:hopf}
Let $X=\mathbb Z_2$ be the biquandle with operations 
$\underline{\triangleright}, \overline{\triangleright}: X \times X \to X$ 
defined by
\[
x\,\underline{\triangleright}\,y= x+1 \mbox{\ \ \  and \ \ \ } x\,\overline{\triangleright}\,y= x+1.
\]
Let $\delta=3 \in \mathbb Z_5$.
Define $A_{xy}, B_{xy}, V_{xy}, C_{xy}, D_{xy}, U_{xy}: X \times X \to \mathbb Z_5$ by
\[
[A_{xy}| B_{xy}| V_{xy}| C_{xy}| D_{xy}| U_{xy}]=
\left[
\begin{array}{cc|cc|cc|cc|cc|cc}
1&1&1&1&0&2&1&3&1&3&0&4\\
1&1&1&1&3&0&3&1&3&1&1&0
\end{array}
\right].
\] 
Note that these are the same biquandle  and biquandle virtual bracket shown in Example~\ref{ex:vtrefoil}.

The right-handed Hopf link $L_1^{\rm Hopf}$ below has nine states with coefficients as listed.
\begin{figure}[htbp]
\begin{center}
  \includegraphics[clip,width=7.0cm]{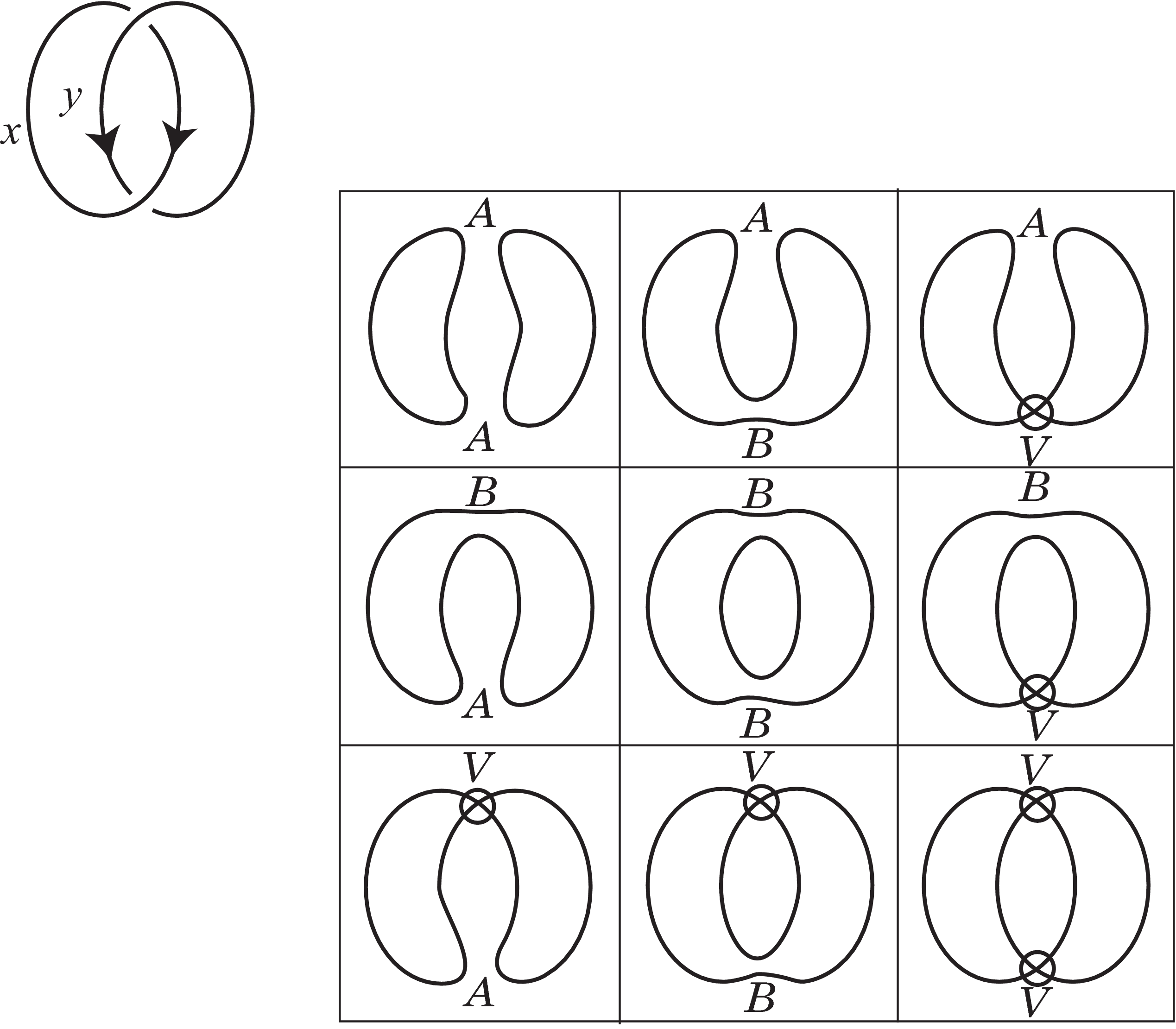}
\end{center}
\end{figure}
Then it has the fundamental biquandle virtual bracket value
\[
\begin{array}{ll}
\phi_1(x,y)=&(A_{xy}A_{yx}+B_{xy}B_{yx}+V_{xy}V_{yx})\delta^2\\
&+(A_{xy}B_{yx}+A_{xy}V_{yx}+B_{xy}A_{yx}+B_{xy}V_{yx}+V_{xy}A_{yx}+V_{xy}B_{yx})\delta.
\end{array}
\]
Hence we have 
\[
\begin{array}{cc|c}
x&y&\phi_1 (x,y) \\ \hline
0&0&(1\cdot1+1\cdot1+0\cdot0)\cdot3^2 +(1\cdot1+1\cdot0+1\cdot1+1\cdot0+0\cdot1+0\cdot1)\cdot3=4\\
0&1&(1\cdot1+1\cdot1+2\cdot3)\cdot3^2 +(1\cdot1+1\cdot3+1\cdot1+1\cdot3+2\cdot1+2\cdot1)\cdot3=3\\
1&0&(1\cdot1+1\cdot1+3\cdot2)\cdot3^2 +(1\cdot1+1\cdot2+1\cdot1+1\cdot2+3\cdot1+3\cdot1)\cdot3=3\\
1&1&(1\cdot1+1\cdot1+0\cdot0)\cdot3^2 +(1\cdot1+1\cdot0+1\cdot1+1\cdot0+0\cdot1+0\cdot1)\cdot3=4
\end{array}
\]
Therefore the right-handed Hopf link $L_1^{\rm Hopf}$ has biquandle virtual bracket invariant
\[
\Phi(L_1^{\rm Hopf})=2 u^{3} + 2 u^{4}.
\]

The left-handed Hopf link $L_2^{\rm Hopf}$ below has nine states with coefficients as listed.
\begin{figure}[htbp]
\begin{center}
  \includegraphics[clip,width=7.0cm]{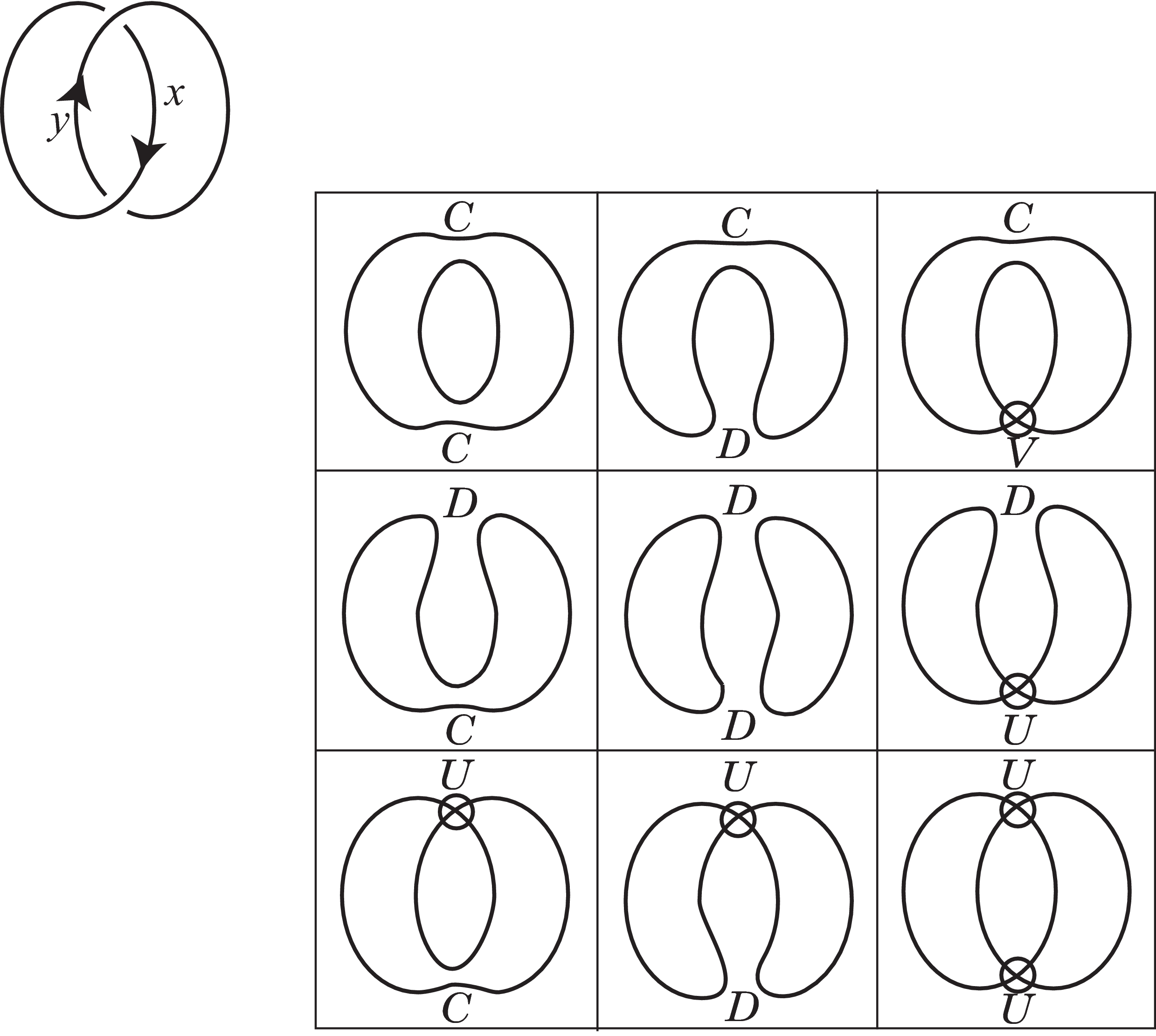}
\end{center}
\end{figure}
Then it has the fundamental biquandle virtual bracket value
\[
\begin{array}{ll}
\phi_2(x,y)=&(C_{xy}C_{yx}+D_{xy}D_{yx}+U_{xy}U_{yx})\delta^2\\
&+(C_{xy}D_{yx}+C_{xy}U_{yx}+D_{xy}C_{yx}+D_{xy}U_{yx}+U_{xy}C_{yx}+U_{xy}D_{yx})\delta.
\end{array}
\]
Hence  we have 
\[
\begin{array}{cc|c}
x&y&\phi_2(x,y)\\ \hline
0&0&(1\cdot1+1\cdot1+0\cdot0)\cdot3^2 +(1\cdot1+1\cdot0+1\cdot1+1\cdot0+0\cdot1+0\cdot1)\cdot3=4\\
0&1&(3\cdot3+3\cdot3+4\cdot1)\cdot3^2 +(3\cdot3+3\cdot1+3\cdot3+3\cdot1+4\cdot3+4\cdot3)\cdot3=2\\
1&0&(3\cdot3+3\cdot3+1\cdot4)\cdot3^2 +(3\cdot3+3\cdot4+3\cdot3+3\cdot4+1\cdot3+1\cdot3)\cdot3=2\\
1&1&(1\cdot1+1\cdot1+0\cdot0)\cdot3^2 +(1\cdot1+1\cdot0+1\cdot1+1\cdot0+0\cdot1+0\cdot1)\cdot3=4
\end{array}
\]
Therefore the left-handed Hopf link $L_2^{\rm Hopf}$ has biquandle bracket invariant
\[
\Phi(L_2^{\rm Hopf})=2 u^{2} + 2 u^{4}.
\]

On the other hand, the unlink of two components $U_2$ has invariant value
\[
\Phi(U_2)=4 u^{4},
\]
and thus, our invariant shows that the right- and left-handed Hopf links are 
non-trivial. 
Moreover, the biquandle virtual bracket invariants of the 
right- and left-handed Hopf links  have the different values, which implies 
that they are not equivalent. In particular, this example shows that some
biquandle virtual bracket invariants are sensitive to orientation reversal.
\end{example}

\begin{example}
Continuing with the same biquandle $X$ from Example \ref{ex:hopf}, let 
$R=\mathbb{F}_8$, the field of eight elements. Recall that $\mathbb{F}_8$
can be written as $\mathbb{Z}_2[t]/\langle 1+t+t^3\rangle=\{0,1,t,1+t,t^2,1+t^2,t+t^2,1+t+t^2\}$. Then our Python computations say that the following defines
a biquandle virtual bracket over $R$:
Let $\delta=1+t$ and 
define $A_{xy}, B_{xy}, V_{xy}, C_{xy}, D_{xy}, U_{xy}: X \times X \to \mathbb{F}_8$ by
\[\begin{array}{l}
[A_{xy}| B_{xy}| V_{xy}| C_{xy}| D_{xy}| U_{xy}]\\
 \\
=\left[
\begin{array}{cc|cc|cc|cc|cc|cc}
1 & 0 & t^2 & 0 & 0 & t^2 & 1 & 0 & 1+t+t^2 & 0 & 0 & 1+t+t^2 \\
0 & 1 & 0 & t^2 & 1+t & 0 & 0 & 1 & 0 & 1+t+t^2 & t+t^2 & 0\\
\end{array}
\right].
\end{array}
\] 
Then our Python computations give the following values for $\Phi_X^{\beta,M}(L)$
for prime classical knots with up to eight crossings:
\[
\begin{array}{l|r}
\Phi_X^{\beta,M}(L) & L \\ \hline
\{2\times 1\} & 5_2, 7_2, 8_{10}, 8_{11}, 8_{13}, 8_{17} \\
\{2\times t\} & 8_3, 8_6, 8_{12}, 8_{116}, 8_{18} \\
\{2\times 1+t\} & U_1, 5_1, 7_6, 8_{15} \\
\{2\times t^2\} & 3_1, 6_2, 8_9 \\
\{2\times 1+t^2\} & 4_1, 7_1, 7_4, 8_5, 8_{14} \\
\{2\times t+t^2\} & 6_1, 6_3, 7_2, 7_3, 8_7, 8_{21} \\
\{2\times 1+t+t^2\} & 7_7, 8_2, 8_3, 8_4, 8_8, 8_{19}, 8_{20},
\end{array}\]
for classical links with up to seven crossings:
\[
\begin{array}{l|r}
\Phi_X^{\beta,M}(L) & L \\ \hline
\{4\times t\} & L2a1 \\
\{4\times 1+t^2\} & U_2 \\
\{2\times 0, 2\times t\} & L7a3 \\
\{2\times 0, 2\times 1+t+t^2\} & L7n2 \\
\{2\times 1, 2\times t^2\} & L6a1 \\
\{2\times t, 2\times 1+t\} & L7a5 \\
\{2\times t, 2\times 1+t^2\} & L7a1, L7a4 \\ 
\{2\times t, 2\times t+t^2\} & L7n1 \\
\{2\times t, 2\times 1+t+t^2 \} & L7a6 \\
\{2\times 1+t, 2\times t^2 \} & L7a2 \\
\{2\times 1+t, 2\times 1+t+t^2\} & L4a1, L6a3, \\ 
\{2\times t^2, 2\times t+t^2 \} & L5a1 \\ 
\{2\times t^2, 2\times 1+t+t^2\} & L6a2 \\
\{8\times t^2\} & U_3 \\
\{2\times 1, 6\times t,\} & L6n1 \\
\{2\times 1+t+t^2, 6\times 1+t^2\} & L6a4 \\
\{2\times t+t^2, 6\times t\} & L6a5 \\
\{2\times t, 6\times 0\} & L7a7 \\
\{2\times t+t^2, 6\times 0\} & L6a5, \\
\{2\times t, 2\times 1+t, 4\times t+t^2\} & L7a7,\\
\end{array}\]
and for virtual knots with up to four classical crossings as listed 
in \cite{KA}:
\[
\begin{array}{l|r}
\Phi_X^{\beta,M}(L) & L \\ \hline
\{2\times 0\} & 3.5, 4.14,4.20,4.21,4.22,4.24,4.34,4.36,4.40,4.52,4.60,4.64, 
4.68,4.89,4.105 \\ \hline
\{2\times t\} & 3.2, 3.3, 3.4, 4.4,4.10,4.11,4.16,4.18,4.23,4.27,4.30,4.31, 4.33,4.38,4.39,4.41, \\ &  4.44,4.45,4.49,4.50,4.54,4.57,4.62,4.63,4.65,4.70,4.74, 4.79,4.81,4.82,4.83,\\ & 4.87,4.92,4.95,4.101 \\ \hline
\{2\times 1+t\} & U_1, 3.1, 3.7, 4.2, 4.6,4.8,4.12,4.13,4.17,4.19,4.26,4.32,4.35,4.42,4.46,4.47,4.51, \\ & 4.55,4.56,4.58,4.59,4.66,4.67,4.71,4.72,4.75,4.76, 4.77,4.85,4.93,4.96,4.97,\\ & 4.98,4.102,4.103,4.106,4.107\\\hline  
\{2\times t^2\} & 3.6 \\\hline
\{2\times 1+t^2\} & 4.1,4.3,4.7,4.9,4.15,4.25,4.29,4.37,4.43,4.48,4.53,4.61, 4.69,4.73,4.78,4.80,\\ & 4.86,4.90,4.91,4.99,4.100,4.108 \\\hline
\{2\times 1+t+t^2\} & 2.1, 4.28,4.84,4.88,4.104.
\end{array}\]

Moreover, $\Phi_X^{\beta}$ distinguishes the square knot $3_1\#\overline{3_1}$ 
from the positive granny knot $3_1\#3_1$, though curiously does not distinguish
the square knot from the
negative granny knot $\overline{3_1}\#\overline{3_1}$ 
\[\includegraphics{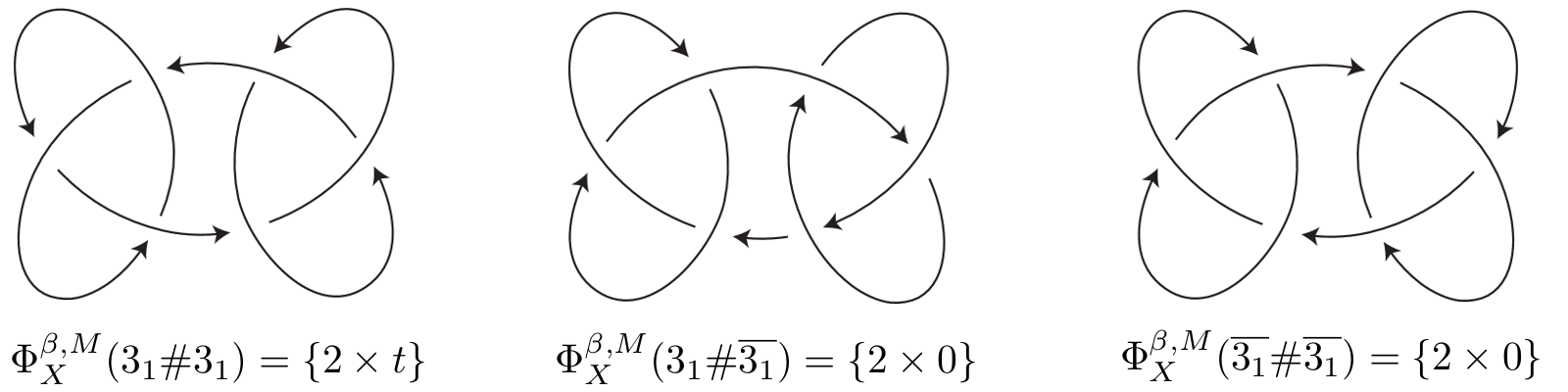}\]
We note that $\Phi_X^{\beta}(\overline{3_1})=\{2\times 0\}$.
\end{example}

\section{\large\textbf{Questions}}\label{Q}

We conclude with some questions for future research.

The case of the square and granny knots example in the last section suggests 
that at least for some virtual brackets, the $\beta$ values might be 
multiplicative
over connected sum; we ask, for which biquandle virtual brackets $(X,R,\beta)$
do we have $\Phi_X^{\beta}(K\#K')=\Phi_X^{\beta}(K)\Phi_X^{\beta}(K')$?

We have only considered the cases of extremely small biquandles and coefficient
rings using computer search; efficient methods for finding virtual brackets 
over larger biquandles and
coefficient rings are of great interest and should produce strong invariants.

As with any invariant, it is worth asking what kinds of categorifications
are possible for these invariants. For Khovanov homology generalizations,
it seems perhaps best to start with biquandle virtual brackets over polynomial
rings.

\bibliography{sn-ko-as-yy-rev}{}
\bibliographystyle{abbrv}

\bigskip

\noindent
\textsc{Department of Mathematical Sciences \\
Claremont McKenna College \\
850 Columbia Ave. \\
Claremont, CA 91711} 

\medskip

\noindent \textsc{Department of Information and Communication Sciences\\
 Sophia University\\
 7-1 Kioi-cho, Chiyoda-ku, Tokyo 102-8554, Japan.
}

\medskip

\noindent \textsc{Department of Mathematics \\
National Institute of Technology, Gunma College \\ 
580 Toriba-cho, Maebashi-shi \\ 
Gunma 371-8530, Japan. }

\end{document}